\RequirePackage[l2tabu, orthodox]{nag}
\documentclass[11pt, reqno]{amsart}

\usepackage{mathrsfs}
\usepackage[dvipsnames]{xcolor}
\usepackage{graphicx}

\usepackage{amsfonts}
\usepackage{amssymb}
\usepackage{amsmath}
\usepackage{amsthm}
\usepackage{amsbsy}
\usepackage{amssymb} 
\usepackage{verbatim}
\usepackage{bm} 
\usepackage{paralist} 
 \usepackage{color}
 \usepackage{mathrsfs}
 \usepackage{graphicx}
\usepackage{hyperref}

\newcommand{\E}{\mathbf{E}}
\def\P{\mathbf{P}}

\def\convd{\stackrel{d}{\rightarrow}}

\def\eqd{\stackrel{d}{=}}

\def\Given{\ \pmb{\big|} \ }

\newcommand{\one}{{\mathbf 1}}

\def\summ{\sum\limits}

\def\to{\rightarrow}

\def\mb{\mbox}

\newcommand{\para}[1]{\vspace{6mm}\noindent{\bfseries #1:}}

\def\l{\left}
\def\r{\right}
\def\<{\langle}
\def\>{\rangle}

\newcommand{\ba}{\[\begin{aligned}}
\newcommand{\ea}{\end{aligned}\]}
\newcommand\mnote[1]{} 
\newcommand{\beq}[1]{\begin{equation}\label{#1}}
\newcommand\eeq{\end{equation}}
\newcommand\ben{\begin{equation}}
\newcommand\een{\end{equation}}
\newcommand\bes{\begin{eqnarray*}}
\newcommand\ees{\end{eqnarray*}}
\newcommand\besn{\begin{eqnarray}}
\newcommand\eesn{\end{eqnarray}}

\def\bthm{\begin{theorem}}
\def\ethm{\end{theorem}}
\def\bdefn{\begin{definition}}
\def\edefn{\end{definition}}
\newcommand{\benu}{\begin{enumerate}\setlength\itemsep{6pt}}
\newcommand{\beit}{\begin{itemize}\setlength\itemsep{3pt}}
\def\eenu{\end{enumerate}}
\def\eeit{\end{itemize}}
\def\beds{\begin{description}}
\def\eeds{\end{description}}
\def\bepr{\begin{problem}}
\def\eepr{\end{problem}}

\def\bprf{\begin{proof}}
\def\eprf{\end{proof}}
\def\berk{\begin{remark}}
\def\eerk{\end{remark}}
\def\bex{\begin{exercise}}
\def\eex{\end{exercise}}
\def\beg{\begin{example}}
\def\eeg{\end{example}}

\newcommand\bllt{$\blacktriangleright$ \;\;}
\def\suchthat{{\; : \;}}
\renewcommand{\qed}{\hfill\text{$\blacksquare$}}

\def\N{\mathbb{N}}
\def\R{\mathbb{R}}

\def\Z{\mathbb{Z}}

\newcommand{\sm}{{\raise0.3ex\hbox{$\scriptstyle \setminus$}}}

\def\alp{\alpha}
\def\bet{\beta}
\def\gam{\gamma}
\def\del{\delta}
\def\eps{\epsilon}
\renewcommand\phi{\varphi}

\def\lam{\lambda}

\def\sig{{\sigma}}

\def\Del{\Delta}

\theoremstyle{plain} 
    \newtheorem{theorem}{Theorem}
    \newtheorem{lemma}[theorem]{Lemma}

\theoremstyle{definition} 
    \newtheorem{definition}[theorem]{Definition}

    \newtheorem{exercise}[theorem]{Exercise}
    \newtheorem{problem}[theorem]{Problem}
        \newtheorem{remark}[theorem]{Remark}
    \newtheorem{example}[theorem]{Example}

\usepackage{charter}

\renewcommand{\benu}{\begin{enumerate}\setlength\itemsep{4pt}}

\def\Sb{\mathbf S}
\def\Vb{\mathbf V}
\def\Ub{\mathbf U}

\def\bas{\small \ba}
\def\eas{\ea \normalsize}

\begin{document}

\title{One idea and two proofs of the KMT theorems}

\author[M. Krishnapur]{Manjunath Krishnapur}
\address{Manjunath Krishnapur, Department of Mathematics, Indian Institute of Science, Bangalore, Karnataka, India.}
\email{manju@iisc.ac.in}

\thanks{The author is partially supported by  UGC Centre for Advanced Study and the SERB-MATRICS grant MTR2017/000292.}

\maketitle


\section{Introduction}

\subsection{The KMT theorems} Koml\'{o}s, Major and Tusn\'{a}dy~\cite{kmt1,kmt2} proved two ``strong embedding'' theorems: one for random walks and   one for empirical processes. 

\bllt {\bf KMT embedding for random walks (KMT-RW)}: Let  $X_{i}$ be i.i.d. random variables with zero mean, unit variance and finite moment generating function in a neighbourhood of zero. Then it is possible to couple the random walk  $S_{k}=X_{1}+\ldots +X_{k}$, with a standard Brownian motion $W$ in such a way that for some constant $C$ and any $x>0$ and any $n\ge 1$, 
\bas
\max_{0\le k\le n}|S_{k}-W(k)| \le C(\log n+x) \;\; \mb{ with probability at least }1-e^{-x}.
\eas

\bllt {\bf KMT embedding for empirical processes (KMT-EP)}: Let $U_{k}$ be i.i.d.\! with $\mb{uniform}[0,1]$ distribution and let $F_{n}(t)=\frac{1}{n}\sum_{k=1}^{n}\one_{U_{k}\le t}$. The {\em uniform empirical process} is the random function $G_{n}(t)=\sqrt{n}(F_{n}(t)-t)$ for $0\le t\le 1$. It is possible to couple a standard Brownian bridge  $W_{0}$  with $G_{n}$ so that for some constant $C$ and any $x>0$ and any $n\ge 1$, 
\bas
\sup_{0\le t\le 1}\sqrt{n}|G_{n}(t)-W_{0}(t)|\le C(\log n+x) \;\; \mb{ with probability at least } 1-e^{-x}.
\eas

These theorems are counted among the most fundamental results of probability theory. There are many extensions to other situations (see the survey by  Lifshits~\cite{lifshits} or the ICM proceedings of Zaitsev~\cite{zaitsev} for some of these) but in this paper we just stick to the versions stated above. The original proofs being rather involved, there have been many efforts to simplify and streamline them. We are aware of two different kinds of proofs. In both, there are two big steps: \begin{inparaenum}[(A)]  \item univariate coupling lemmas and \item extension of  the coupling to the level of paths. \end{inparaenum} We describe the two approaches next. To see the results, one may jump directly to Sections~\ref{sec:univariatecouplinglemmas} and\ref{sec:ourapproach}.

\subsection{The  two approaches} The first one is the original proof of Koml\'{o}s, Major and Tusn\'{a}dy with further developments due to Cs\"{o}rg\~o-R\'{e}v\'{e}sz~\cite{csorgorevesz}, Bretagnolle and Massart, Dudley, Massart (see \cite{massart} for a discussion of these papers), Carter-Pollard~\cite{carterpollard} and Pollard~\cite{pollard} (this is far from a complete list, see references in \cite{pollard} and \cite{carterpollard}). This   works for both versions of KMT.

Step~(A) in  KMT-RW consists in proving a version of Tusn\'{a}dy's lemma (Lemma~\ref{lem:snandzcoupling}), which gives a coupling of $S_{n}$ with a Gaussian. For KMT-EP, one needs this lemma only for Bernoulli steps (the connection is that $nF_{n}(t)\sim \mb{Binomial}(n,t)$). This is achieved by a fine comparison of the tail probabilities of $S_{n}$ with that of a Gaussian, using Stirling's approximation.

Step~(B): For KMT-RW, the proof uses Step~(A) to couple $S_{n}$ with $W(n)$, then $S_{n/2}$ with $W(n/2)$, then $S_{n/4}$ and $S_{3n/4}$, etc. For KMT-EP, the proof couples $G_{n}(t)$ and $W_{0}(t)$ for dyadic $t$. In either case, the dyadic procedure is carried up to a depth of about $\log n$ generations. An excellent exposition of the complete proof of the KMT theorem for empirical processes is the book by Pollard~\cite{pollard}.

 The second approach, due to Chatterjee~\cite{chatterjee}, proves KMT-RW when the steps have symmetric Bernoulli distribution, by an approach that may be broadly described as Stein's method. Bhattacharjee and Goldstein~\cite{bhattacharjeegoldstein} extended this method of proof to a large class of step distributions.

Step~(A) in this method consists in constructing couplings of Binomial and Hypergeometric distributions with Gaussian distributions.  This is achieved by constructing the Stein coefficients for these distributions (Stein coefficient is a tool that measures how far a distribution or a random variable is from satisfying Stein's equation for the Gaussian), and showing that they are close to constants. A key step is a general new result obtained by Chatterjee that deduces from this the existence of a good coupling. These couplings are weaker than Tusn\'{a}dy's lemma, but suffice for the next step.

Step~(B) comprises of a  non-trivial induction on the number of steps of the random walk. To make the induction work, the hypothesis chosen is a careful statement about coupling the random walk bridge  with a Brownian bridge having the same endpoint. 

 The goal of this paper is to give new and possibly simpler ways to carry out Step~(A) in both methods, for the symmetric Bernoulli case. We have little to add to Step~(B),  but for completeness, we sketch them in the two appendices. The appendices are not original and are taken from  Chapter 10 of Pollard's book~\cite{pollard} for KMT-EP and Sections~4,5 of Chatterjee's paper~\cite{chatterjee} for the first KMT-RW (symmetric Bernoulli case only). 

\subsection{The univariate coupling lemmas}\label{sec:univariatecouplinglemmas} We now state the three Lemmas that we prove in this paper. These comprise Step~(A) in the two proofs outlined above. Throughout the paper, unless indicated otherwise, $Z$ will denote a standard Gaussian random variable and $S_{n}=X_{1}+\ldots +X_{n}$ where $X_{i}$ are i.i.d. symmetric Bernoullis, i.e.,  $X_{i}=\pm 1$ with probability $1/2$ each.

\subsubsection{\textbf{Step~(A) in the proof of KMT-EP}} This  consists entirely in the following famous lemma of Tusn\'{a}dy (the usual reference is to \cite{tusnadythesis} but what we  understand comes from Pollard~\cite{pollard}). 
\begin{lemma}[Tusn\'{a}dy-type lemma]\label{lem:snandzcoupling} For some $n_{0}$ and any $n\ge n_{0}$, there is a coupling of $S_{n}$ with $Z$ such that $|S_{n}|\le |Z|\sqrt{n}+3$ and  $|S_{n}-Z\sqrt{n}|\le Z^{2}+11$.
\end{lemma}
In the strongest known form of Tusn\'{a}dy's lemma of this form, the result is valid for all $n\ge 1$ with $|Z|\sqrt{n}+2$ and $\frac14 Z^{2}+2$ on the the right sides of the two inequalities. Carter and Pollard~\cite{carterpollard} improved the main term on the right side of the second inequality to $C_{1}+C_{2}Z^{2}(1\wedge \frac{1}{n}|Z|)$, which is smaller for typical values of $Z$. 

Our version has explicit and decent constants, and although suboptimal, they can be brought down further, see Section~\ref{sec:bestconstant}. We are not aware of any use of these constants though. In fact, in deriving the KMT for empirical processes, Pollard~\cite{pollard}  (chapter~10 and Appendix~D) uses a weaker lemma in which the right sides are $C(1+|Z|\sqrt{n})$ and $C(1+Z^{2})$ in the two inequalities, with unspecified constants. Of course,  the statement of the Lemma~\ref{lem:snandzcoupling} implies the same for all $n\ge 1$, albeit with  larger constants.

\subsubsection{\textbf{Step~(A) in the proof of KMT-RW for symmetric Bernoulli steps}}
Chatterjee proves the following lemmas in place of Tusn\'{a}dy's and uses them to derive the KMT theorem for the Bernoulli random walk. These are Theorem~3.1 and Theorem~3.2 in \cite{chatterjee}. 

The first lemma is about coupling a binomial distribution with a Gaussian distribution of the same mean and variance. Like in Tusn\'{a}dy's lemma, the distance between the two does not grow with the variance of the variables.
\begin{lemma}[Theorem~3.1 in \cite{chatterjee}]\label{lem:bingausscoupling} For some $\theta_{0}>0$ and $\kappa_{0}<\infty$, for  every $n\ge 1$, there exists a coupling of $S_{n}$ with $Z$ such that $\E[e^{\theta_{0}  |S_{n}-Z\sqrt{n}|}]\le \kappa_{0}$. 
\end{lemma}
In fact, our proof allows any $\theta_{0}$ satisfying $8e^{2\theta_{0}}\theta_{0}^{2}<1$ and gives an explicit form for $\kappa_{0}$.

The second  lemma that Chatterjee proves is a coupling of a hypergeometric distribution with a Gaussian distribution. The relevance of this is easy to see: In Step~(B), when one proceeds by conditioning on position of the random walk at the end, the position of the walk at any intermediate time has a (shifted) hypergeometric distribution. 

Some notation: For $0\le k\le n$ and any probable value $s$ of $S_{n}$ (by that we mean $\P\{S_{n}=s\}>0$, or equivalently, that $|s|\le n$ and $n-s$ is even),  let  $S_k[n,s]$ denote a random variable whose distribution is the conditional distribution of $S_{k}$ given $S_{n}=s$. An equivalent description is that $S_k[n,s]$ is the sum of the first $k$ coupons drawn without replacement, uniformly at random, from a box containing $n$ coupons of which $p=(n+s)/2n$ proportion are labeled $+1$ and the remaining $q=(n-s)/2n$ proportion are labeled $-1$.  Thus $S_k[n,s]$ has  (a the simple transformation of)  hypergeometric distribution and has mean $k(p-q)=\frac{sk}{n}$ and variance $4pq\frac{k(n-k)}{n-1}$.  Let $\sig_{n,k}^{2}=\frac{1}{n}k(n-k)$.

\begin{lemma} \label{lem:hypergeogausscoupling} There exists a $\theta_{1}>0$ and $M_{1}<\infty$ such that for any   $n\ge 2$ and $\frac13 n \le k\le \frac23 n$, and any probable value $s$ of $S_{n}$, there exists a coupling of $W:=S_k[n,s]-\frac{sk}{n}$ with $Z\sim N(0,1)$ such that $\E[e^{\theta|W-\sig_{n,k} Z|}]\le \exp\{1+M_{1}\theta^{2}\frac{s^{2}}{n}\}$ for all $\theta\le \theta_{1}$.
\end{lemma}
Note that $\sig_{n,k}^{2}$ does not depend on $s$ and  is comparable to the variance of $S_{k}[n,s]$ only when $s=O(\sqrt{n})$. For such $s$, the conclusion here is analogous to that of Lemma~\ref {lem:bingausscoupling}. For larger values of $s$, the  coupling here  is not between random variables of comparable variance and correspondingly, the right side could be large, but the point is that the controlling parameter is $s^{2}/n$.  

\subsection{Our approach}\label{sec:ourapproach} The single key idea of this paper is that {\em Binomial distributions are similar objects (discrete, combinatorial), and comparing them with each other is easier than comparing one of them to the Gaussian distribution}.  This suggests that we approach the problem via the {\em Cauchy criterion} and look for a coupling between $S_n/\sqrt{n}$ and $S_{4n}/\sqrt{4n}$.  Then we can successively couple\footnote{When we talk about a coupling of random variables $X,Y$ satisfying some statements, it means that there exists some probability space on which there are random variables having the same marginal distributions as $X$ and $Y$, and satisfying the mutual relationships in the statements.  Given couplings $\theta_{i}$ of $X_{i}$ and $X_{i+1}$ for $i=0,1,\ldots$, one can get a sequence of random variables $(X_{i}')_{i\ge 0}$ such that $(X_{i}',X_{i+1}')$ has the joint distribution given by the coupling $\theta_{i}$ for each $i\ge 0$. This is done by running a Markov chain as follows: First sample $(X_{0}',X_{1}')$ from the distribution $\theta_{0}$, and then successively sample $X_{i+1}'$ from the conditional distribution of $X_{i+1}$ given $X_{i}$ in the distribution $\theta_{i}$.
} $S_{n}/\sqrt{n}$, $S_{4n}/\sqrt{4n}$, \ldots ,$S_{n4^k}/\sqrt{n4^k},\ldots$.  If the coupling is sufficiently strong at each step, this sequence converges almost surely  to a standard Gaussian variable $Z$ that is  coupled well with $S_n/\sqrt{n}$. 

This is clearly an approach that can be of more general use, for instance to get rates of convergence in limit theorems where often discrete combinatorial objects converge to a continuum object. In our context, this idea translates to proving the following theorems, which are natural analogues of Lemmas~\ref{lem:snandzcoupling}, \ref{lem:bingausscoupling}, \ref{lem:hypergeogausscoupling}.
\begin{theorem} \label{thm:snands4ncoupling} For some $n_{0}$ and any even number $n\ge n_{0}$, there exists a coupling of $2S_{n}$ and $S_{4n}$ so that they have the same sign (meaning $S_{n}S_{4n}\ge 0$) and 
\bas
|S_{4n}|-\frac{1}{8n}|S_{4n}|^{2}-9 \; \le \;  2|S_{n}| \; \le \; |S_{4n}|+2.
\eas
\end{theorem} 
Equivalently,  $2|S_{n}|\le |S_{4n}|+2$ and $|2S_{n}- S_{4n}|\le \frac{1}{8n}|S_{4n}|^{2}+9$. This looks more  similar to Lemma~\ref{lem:snandzcoupling}, with $\frac{1}{\sqrt{4n}}S_{4n}$ taking the place of $Z$. 
\begin{theorem}\label{thm:bingausscouplingcauchy} There exists $\theta_{0}>0$ and $\kappa_{0}<\infty$ such that  for any  $n\ge 1$, there exists a coupling of $S_{n}$ with $S_{4n}$ such that $\E[e^{\theta_{0} |2S_{n}-S_{4n}|}]\le \kappa_{0}$. In fact, any $\theta_{0}$ such that $8\theta_{0}^{2}e^{2\theta_{0}}<1$ works.
\end{theorem}

The analogue of Lemma~\ref{lem:hypergeogausscoupling} is broken into two parts. In the unbiased case $s=0$, we show that $S_k[n,0]$ can be coupled well with a Gaussian of the same variance. Then we show that $S_k[n,s]$ can be coupled with $S_k[n,0]$ so that the difference is controlled by $s^{2}/n$.
\begin{theorem}\label{thm:hypergeogausscouplingcauchy}  There exists a $\Theta>0$ and $M<\infty$ such that for any  even number $n$ and $\frac13 n \le k\le \frac23 n$,  and any probable value $s$ of $S_{n}$, writing $W_{1}:=S_k[n,0]$, $W_{2}:=S_{4k}[4n,0]$ and $W:=S_k[n,s]-\frac{sk}{n}$, there are couplings
\benu
\item   of $W_{1}$ with $W_{2}$ such that $\E[e^{\theta|2W_{1}-W_{2}|}]\le \frac32$ for all $\theta\le \Theta$, and
\item  of $W_{1}$ with $W$ such that $\E[e^{\theta|W_{1}-W|}]\le e^{1+M\theta^{2}\frac{s^{2}}{n}}$ for all $\theta\le \Theta$.
\eenu
\end{theorem}
We assume $n$ to be even so that $0$ is a probable value of $S_{n}$ and hence $S_k[n,0]$ makes sense. For odd $n$, the same holds if we repalce $0$ by $1$ (or any fixed odd number). In another direction, we could have stated the theorem more generally for  $\del n\le k\le (1-\del)n$ for any $\del>0$. Alternately from the point of view of its actual use in Step~(B), we could have specialized to $k=\lfloor n/2\rfloor$ and simplified the proofs a little.  Since Chatterjee wrote it for $\frac{k}{n}\in [\frac13,\frac23]$, we do the same.

\subsection{About the proofs}
The potential to use combinatorial methods eliminates many technicalities to make the proofs of Theorems~\ref{thm:snands4ncoupling}, \ref{thm:bingausscouplingcauchy}, \ref{thm:hypergeogausscouplingcauchy}  simpler relative to Lemmas~\ref{lem:snandzcoupling}, \ref{lem:bingausscoupling}, \ref{lem:hypergeogausscoupling}. And the deduction of those three lemmas from the corresponding theorems is also straightforward. While the  simplicity may be in the eye of the beholder, given the fundamental nature of the KMT theorems, even improvement of exposition may be of some interest. Carter and Pollard~\cite{carterpollard} remark on the ``...continuing perceived need for an accessible treatment of  the coupling result that underlies the KMT construction''. 

\subsubsection{About the proof of Theorem~\ref{thm:snands4ncoupling}} The main step are Lemmas~\ref{lem:betbeatsalp} and \ref{lem:alpbeatsbet}, which give a comparison of Binomial coefficients of the form $\binom{n}{i}$ with $\binom{4n}{j}$. This  comparison of the Binomial mass functions, together with some standard estimates on the tails of Binomial distribution allows us to prove the right kinds of inequalities between the tails of the two Binomial distributions. In contrast, earlier proofs (as in Pollard~\cite{pollard}, Carter and Pollard~\cite{carterpollard}, Massart~\cite{massart}, etc.) directly compare the Binomial tails to Gaussian tails. The estimates needed are more refined than the usual first-order Stirling's formula used to prove the de Moivre-Laplace central limit theorem. 

\subsubsection{About the proofs of Lemma~\ref{lem:bingausscoupling} and Lemma~\ref{lem:hypergeogausscoupling}} First we outline Chatterjee's approach - all references are to \cite{chatterjee}. If there is a random variable $T$ (usually $T=T(W)$, although use of additional randomness is allowed) such that $\E[Wf(W)]=\E[Tf'(W)]$ for a large class of $f$, then Chatterjee calls $T$ a Stein coefficient for $W$. Stein's famous characterization of the Gaussian says that if $T=\sig^{2}$ is a constant, then $W\sim N(0,\sig^{2})$. A key result in Chatterjee's paper is that if $T$ is close to a constant in an appropriate sense, then $W$ is close to a Gaussian. Of course, this statement needs to be proved in a much finer form than usually needed to show central limit theorems by Stein's method. After that, it is a matter of finding the  Stein coefficients of the Binomial and Hypergeometric distributions (actually of some perturbations of those, since they are discrete and do not admit Stein coefficients) to prove  Lemma~\ref{lem:bingausscoupling} and Lemma~\ref{lem:hypergeogausscoupling}.

In our proofs of Theorem~\ref{thm:bingausscouplingcauchy} and Theorem~\ref{thm:hypergeogausscouplingcauchy}, we work with finite state Markov chains and show that if two  Markov chains on (segments of) integers have transitions $i\mapsto i\pm 1$ with rates $T(i)\mp i$ and $i\mapsto i\pm 1$ with rates $S(i)\mp i$, then if $T$ and $S$ are close, the stationary distributions of the two chains can be coupled well.  This is achieved by constructing a joint Markov chain on $\Z^{2}$, whose stationary distribution gives a coupling of the stationary distributions of the given chains.

Although this may seem different from what was outlined as Chatterjee's method, they are closely related. Stein's equation is essentially a rephrasing of the forward equation for the Ornstein-Uhlenbeck process. And the proof of the main lemma of Chatterjee essentially boils down to  showing the existence of an invariant measure for a certain generator of a joint Markov process. The technicalities are considerably reduced in our case as we are in the setting of finite state space Markov chains.

An alternate way to say the same thing is that we work with the Stein operator for the Binomial  distributions. It is a fact that 
\bas
\E\l[\l(\frac{n}{2}+X\r)(f(X-1)-f(X))\r]=\E\l[\l(\frac{n}{2}-X\r)(f(X+1)-f(X))\r]
\eas
for all functions $f$ if and only if $X$ has the same distribution as $S_{n}$. Thus, to prove Theorem~\ref{thm:bingausscouplingcauchy}, we show that $Y=2(S_{n}+R)$ (where $R$ is independent of $S_{n}$ and takes values $-1,0,1$) satisfies
\bas
\E[(T(Y)+Y)(f(Y-1)-f(Y))]=\E[(T(Y)-Y)(f(Y+1)-f(Y))]
\eas
for all $f$, for a function $T$ that is close to the constant $2n$. This shows that $Y$ (and hence $2S_{n}$)  can be coupled well with $S_{4n}$.

\subsubsection{A minor point worth noting} In his famous paper where the well-known inequality comes, Hoeffding~\cite{hoeffding} also proved other less known but remarkable results. One of them states that if $f:\R\mapsto \R$ is convex, then $\E[f(\bar{X})]\le \E[f(\bar{Y})]$ where $\bar{X}$ and $\bar{Y}$ are averages of samples drawn without and with replacement from a box of coupons. This allows one to get estimates on expectations of functions of a hypergeometric variable (see Lemma~\ref{lem:boundsinbinomialandhypergeo}) without much effort. Not using this result of Hoeffding, direct and lengthier proofs are given in  \cite{chatterjee} (Lemmas~3.4 and 3.5 in \cite{chatterjee} and Lemma~\ref{lem:boundsinbinomialandhypergeo} in this paper).

\subsubsection{Some shortcomings} The first proof, although mainly combinatorial, could be made nicer if one could prove the main Lemmas~\ref{lem:betbeatsalp} and \ref{lem:alpbeatsbet} by bijective methods. We were unable to do that. Further, Lemma~\ref{lem:alpbeatsbet} strongly suggests that it should be possible to prove the imroved version of Tusn\'{a}dy's lemma as found by Carter and Pollard. However, in trying to do that the remaining part of the proof got so bloated that we settled for the weaker form.

The second proof, by Markov chain coupling, appears to be amenable to proving the KMT theorem for random walks with more general step distributions. The reason is that the Stein coefficient has nice behaviour under convolutions, and it appears that Cramer's large deviation theorem should give a coupling of $S_{n}$ with $N(0,n)$ analogous to Theorem~\ref{thm:bingausscouplingcauchy}. Of course, this is only one of the key steps and there is more work needed.  As of now, we do not have such a proof. As was mentioned earlier, Bhattacharjee and Goldstein~\cite{bhattacharjeegoldstein} have already extended Chatterjee's method to more general random variables, although with some extra conditions.

\subsection{Outline of the rest of the paper}  In the next section, we deduce Lemmas~\ref{lem:snandzcoupling}, \ref{lem:bingausscoupling}, \ref{lem:hypergeogausscoupling} from Theorems~\ref{thm:snands4ncoupling}, \ref{thm:bingausscouplingcauchy} and \ref{thm:hypergeogausscouplingcauchy}, respectively. Then the paper is split into two parts that can be read independently of each other. Part~\ref{partI} proves Theorem~\ref{thm:snands4ncoupling} by a combinatorial method and Part~\ref{partII} proves Theorems~\ref{thm:bingausscouplingcauchy} and \ref{thm:hypergeogausscouplingcauchy} by coupling Markov chains. The overall ideas of the proofs were outlined above. The appendices outline Step~(B) in both proofs, i.e., the deduction of the KMT theorems from the  coupling lemmas.

\subsection{Acknowledgments} I first learned about the KMT theorem from Yuval Peres when (around 2007)  he was asking for simpler proofs. That motivated me to think about the problem off and on. I would like to thank Yogeshwaran D. for listening  to my speculative ideas at various stages, for useful comments on the first draft, and for encouraging remarks at all times.

\section{Proofs of Lemmas~\ref{lem:snandzcoupling}, \ref{lem:bingausscoupling}, \ref{lem:hypergeogausscoupling} from Theorems~\ref{thm:snands4ncoupling}, \ref{thm:bingausscouplingcauchy} and \ref{thm:hypergeogausscouplingcauchy}}
All three proofs are similar. For Lemmas~\ref{lem:snandzcoupling}, \ref{lem:bingausscoupling}, we define $Z_{j}=S_{4^{j}n}/\sqrt{4^{j}n}$, for $j=0,1,2\ldots$. From the corresponding theorems, we can  couple $Z_{j}$ with $Z_{j+1}$ for each $j$, and hence a couple all the $Z_{j}$s on a common probability space so that $Z_{j}$ and $Z_{j+1}$ are close for each $j$ (with high probability). Then we show that the sequence $Z_{j}$ converges almost surely to a standard Gaussian random variable $Z$, and that $Z_{0}$ and $Z$ are very close. That is $S_{n}$ and $Z$ are coupled as required by the Lemmas. For Lemma~\ref{lem:hypergeogausscoupling}, the idea is quite similar but the notation is different and will be introduced in the proof.

\bprf[Proof of Lemma~\ref{lem:snandzcoupling}] Fix an even number $n\ge n_{0}$ as in Theorem~\ref{thm:snands4ncoupling}, and couple all the   $Z_{k}$s  so that  they all have the same sign and 
\small
\begin{align}\label{eq:couplingofZks} 
|Z_{k}|\le |Z_{k+1}|+\frac{1}{2^{k}\sqrt{n}}, \qquad |Z_{k}-Z_{k+1}|\le \frac{1}{2^{k+2}\sqrt{n}}|Z_{k+1}|^{2}+\frac{9}{2^{k+1}\sqrt{n}}.
\end{align}
\normalsize
By Bernstein's inequality, $\P\{|Z_{k+1}|\ge 2^{k/4}\}$ is summable, hence $|Z_{k+1}|\le 2^{k/4}$ for all but finitely many $k$. Then the second inequality in \eqref{eq:couplingofZks} shows that $|Z_{k}-Z_{k+1}|$ is summable, and hence $(Z_{k})_{k}$ is a Cauchy sequence, almost surely. The limiting variable, call it $Z$, must have $N(0,1)$ distribution, by the central limit theorem. Evidently, $Z$ has the same sign as $Z_{0}$.

Summing the first inequality in \eqref{eq:couplingofZks} over $k\ge j$ gives $|Z_{j}|\le |Z|+2^{1-j}n^{-1/2}$. In particular  $|Z_{0}|\le |Z|+2n^{-1/2}$ or equivalently $|S_{n}|\le |Z|\sqrt{n}+2$. 

Use $|Z_{k+1}| \le |Z|+2^{-k}n^{-\frac12}$ in the second inequality in \eqref{eq:couplingofZks} and sum  to get
\small
\begin{align}
|Z_{0}-Z| &\le \sum_{k=0}^{\infty}\frac{1}{2^{k+2}\sqrt{n}}\l(|Z|+\frac{1}{2^{k}\sqrt{n}}\r)^{2}+\frac{9}{2^{k+1}\sqrt{n}} \nonumber \\
&\le \frac{1}{\sqrt{n}}Z^{2}+\frac{4}{7n^{3/2}} + \frac{9}{\sqrt{n}}  \label{eq:wastefulboundonsquareofsum}
\end{align}
\normalsize
where we used $(a+b)^{2}\le 2a^{2}+2b^{2}$ and summed the geometric series. Consequently,  $|Z_{0}-Z|\sqrt{n}\le Z^{2}+10$ for large $n$. Equivalently, $|S_{n}-Z\sqrt{n}|\le Z^{2}+10$. 

This completes the proof of the Lemma for even $n$, with a saving of $1$ on the right sides of both inequalities.  If $n$ is odd, then $S_{n}$ can be coupled with $S_{n+1}$ (by adding one step in the random walk) so that $|S_{n}-S_{n+1}|\le 1$, and using the coupling of $S_{n+1}$ with $Z$, this proves the Lemma for  odd $n$.
\eprf
 
Next we deduce Lemma~\ref{lem:bingausscoupling} from Theorem~\ref{thm:bingausscouplingcauchy}. 
\bprf[Proof of Lemma~\ref{lem:bingausscoupling}] Fix any $\theta_{0}$ and $\kappa_{0}$ as in Theorem~\ref{thm:bingausscouplingcauchy}. Couple the variables $Z_{0},Z_{1},Z_{2},\ldots$ on the same probability space so that for each $k\ge 0$
\bas
\E\l[e^{\theta_{0}\sqrt{n}2^{k+1}|Z_{k}-Z_{k+1}|}\r]\le \kappa_{0}.
\eas
By convexity of the exponential,
\bas
\exp\{\theta_{0}\sqrt{n} \sum_{k\ge 0}|Z_{k}-Z_{k+1}|\} \le \sum_{k\ge 0}\frac{1}{2^{k+1}}\exp\l\{\theta_{0}\sqrt{n}2^{k+1}|Z_{k}-Z_{k+1}|\r\}.
\eas
Therefore, the expectation of the left side quantity is bounded by $\kappa_{0}$. In particular, this shows that $\sum_{k\ge 0}|Z_{k}-Z_{k+1}|$ converges almost surely, which  implies that $Z_{k}$ converges almost surely to some $Z$. By the central limit theorem, $Z\sim N(0,1)$. As $|Z_{0}-Z|\le \sum_{k\ge 0}|Z_{k}-Z_{k+1}|$ it follows that $\E[e^{\theta_{0} \sqrt{n}|Z_{0}-Z|}]\le \kappa_{0}$. This is  the same as $\E[e^{\theta_{0}|S_{n}-Z\sqrt{n}|}]\le \kappa_{0}$. 
\eprf

\bprf[Proof of Lemma~\ref{lem:hypergeogausscoupling}]  Recall that $S_k[n,0]$ has mean $0$ and variance $\frac{n}{n-1}\sig_{n,k}^{2}$. Set $U_{j}=2^{-j}S_{4^{j}k}[4^{j}n,0]$. Then $U_{j}$ has zero mean and  variance $(1-\frac{1}{4^{j}n})^{-1}\sig_{n,k}^{2}$.  From the first part of Theorem~\ref{thm:hypergeogausscouplingcauchy}, if $n$ is an even number, we can construct random variables $U_{j}$ on one probability space so that for any $j\ge 0$ and any $\theta\le \Theta$, we have
$\E[e^{\theta 2^{j+1}|U_{j}-U_{j+1}|}]\le \frac32$.  Using convexity of the exponential, 
\bas
\E\l[e^{\theta \sum_{j\ge 0}|U_{j}-U_{j+1}|}\r] &\le \E\l[\sum_{j\ge 0}\frac{1}{2^{j+1}}e^{\theta 2^{j+1}|U_{j}-U_{j+1}|}\r] \; \le \frac32
\eas
for $\theta\le \Theta$. This shows that $\sum_{j}|U_{j}-U_{j+1}|$ converges almost surely and hence $U_{j}$ converges almost surely to some random variable $U$, as $j\to \infty$. As $|U_{0}-U|$ is bounded by $\sum_{j\ge 0}|U_{j}-U_{j+1}|$, we see that $\E[e^{\theta|U_{0}-U|}]\le \frac32$ for $\theta\le \Theta$.  Further, this also shows that $\{U_{j}\}$ is exponentially tight and hence its mean and variance converge to those of $U$. Since $U_{j}$ converges in distribution to $N(0,\sig_{n,k}^{2})$  (this is an elementary fact, see Remark~\ref{rem:cltforsrswor}), it follows  that $U\sim N(0,\sig_{n,k}^{2})$. In other symbols, we have a coupling of $W_{1}=U_{0}$ with $Z:=U/\sig_{n,k}$ such that $\E[e^{\theta|W_{1}-\sig_{n,k}Z|}]\le \frac32$ for all $\theta\le \Theta$.

Next, use the second part of Theorem~\ref{thm:hypergeogausscouplingcauchy} to construct $W=S_{k}[n,s]-\frac{sk}{n}$ coupled with $W_{1}$ in such a way that $\E[e^{\theta|W-W_{1}|}]\le \exp\{1+M\theta^{2}\frac{s^{2}}{n}\}$ for all $\theta\le \Theta$. 

Now that we have $W,W_{1},Z$ on the same probability space, we just observe that for $\theta\le \frac12 \Theta$, by Cauchy-Schwarz,
\bas
\E[e^{\theta|W-\sig_{n,k}Z|}] 
&\le \E[e^{2\theta|W_{1}-\sig_{n,k}Z|}]^{\frac12}\ \E[e^ {2\theta|W-W_{1}|}]^{\frac12} \\
&\le \sqrt{\frac32} \ e^{\frac12+2M\theta^{2}\frac{s^{2}}{n}} \\
\eas
As this is less than $e^{1+2M\theta^{2}\frac{s^{2}}{n}}$,  taking $\theta_{1}=\Theta/2$ and $M_{1}=2M$, this completes the proof for even $n$.

If $n$ is odd, we claim that $S_{k}[n,s]$ can be coupled with $S_{k}[n+1,s+\frac12]$ so that the difference between the two is at most $2$. To see this, consider the box of $(n+s)/2$ coupons labeled $1$ and $(n-s)/2$ coupons labelled $-1$ and a red coupon also labelled $+1$. Drawing $k$ coupons without replacement and adding them gives $S_{k+1}[n,s]$. To get $S_{k}[n,s]$, we do the same experiment, but if the red coupon comes up, discard it and draw a different one. It is clear that this can be done so that the difference between the two remains at most $2$. Hence, using the result for even $n$, we get a coupling such that
\ba
\E[e^{\theta|W-\sig_{n,k}Z|}] &\le e^{2\theta}\sqrt{\frac32} \ e^{\frac12+2M\theta^{2}\frac{(s+\frac12)^{2}}{n}}
\ea 
which is bounded by $e^{1+4M\theta^{2}\frac{s^{2}}{n}}$ if $\theta$ is small enough. 
\eprf
\berk\label{rem:cltforsrswor} We claimed that $U_{j}\convd N(0,\sig_{n,k}^{2})$ in the proof. This follows from the fact that if $t\in (0,1)$ is  fixed, $V_{n}:=S_{\lfloor nt\rfloor }[n,0]/\sqrt{n}$ converges in distribution to $N(0,t(1-t)/4)$. A direct way to show this is to use Stirling's approximation to write the probability that  $V_{n}=x$ as (here $H(t)=-t\log t -(1-t)\log(1-t)$)
\bas
& \frac{\binom{m}{mt+\frac12x\sqrt{2m}}\binom{m}{mt-\frac12x\sqrt{2m}}}{\binom{2m}{2}} &\sim \frac{1}{\sqrt{\pi m}}e^{m[H(t+\frac{x\sqrt{2}}{\sqrt{2m}})+H(t-\frac{x\sqrt{2}}{\sqrt{2m}})-2H(t)]}.
\eas
The exponent is $-\frac{2x^{2}}{t(1-t)}+o(1)$, as the terms linear in $x$ cancel. Thus, $V_{n}$ converges in distribution $N(0,\frac14t(1-t))$.
\eerk

\newpage
\part{Coupling by comparison of binomial probabilities}\label{partI}

\section{The comparison between mass functions of $2S_{2m}$ and $S_{8m}$}
In this section, assume that $n=2m$ is even. Then $2S_{2m}$ is supported on multiples of $4$ from $-2m$ to $2m$ and $S_{8m}$ is supported on even integers from $-8m$ to $8m$. Set $\alp_{m}(k)=\P\{2S_{2m}=4k\}$ and   $\bet_{m}(k)=\P\{S_{8m}=4k\mb{ or } 4k-2\}$ for $k\ge 1$. Explicitly,
\bas
\alp_{m}(k) &=\binom{2m}{m+k}\frac{1}{2^{2m}} \\
\bet_{m}(k) &=\binom{8m}{4m+2k}\frac{1}{2^{8m}}+\binom{8m}{4m+2k-1}\frac{1}{2^{8m}} = \binom{8m+1}{4m+2k}\frac{1}{2^{8m}}.
\eas
 This are not mass functions, but if we  double the values of $\alp_{m}(k)$ and $\bet_{m}(k)$ for $k\ge 1$ and set $\alp_{m}(0)=\P\{S_{2m}=0\}$ and $\bet_{m}(0)=\P\{S_{8m}=0\}$, we get the mass functions of $2|S_{2m}|/4$ and $\lceil |S_{8m}|/4\rceil$. The ultimate goal is to couple $2S_{2m}$ and $S_{8m}$ as closely as possible, which is possible if one gets a good comparison between the tails of the two distributions (see Lemma~\ref {lem:gencoupling}). As a first step towards this goal, in this section we get a comparison between the probability mass functions.
\begin{lemma}\label{lem:betbeatsalp} For any $m\ge 1$, we have $\alp_{m}(k)\le \bet_{m}(k)$ for any $k\ge 1$.
\end{lemma}
\begin{lemma}\label{lem:alpbeatsbet} For any $m\ge 1$, we have $\alp_{m}(k)\ge \bet_{m}(\ell)$ if $1\le k\le \ell-\frac14\l(1+\frac{\ell^{3}}{m^{2}}\r)$. In particular, $\alp_{m}(\ell-1)\ge \bet_{m}(\ell)$ for $\ell\le (3m^{2})^{1/3}$.
\end{lemma}
We now proceed to the proofs of the two lemmas, which are similar. For $h\ge 1$, define
\bas
f(m,k)&:=\frac{\bet_{m}(k)}{\alp_m(k)} =\frac{\binom{8m+1}{4m+2k}}{2^{6m}\binom{2m}{m+k}}, \\ 
g_{h}(m,k)&:=\frac{\bet_{m}(k)}{\alp_{m}(k-h)} \; = \;\frac{\binom{8m+1}{4m+2k}}{2^{6m}\binom{2m}{m+k-h}}.
\eas
The assertion of the two lemmas  is that $f(m,k)\ge 1$ for $1\le k\le m$ and $g_{h}(m,k)\le 1$ if $h+1\le k\le m$ and $h\ge \frac14(1+\frac{k^{3}}{m^{2}})$.  It is tempting to try proving these by finding explicit injective maps between appropriate sets, since the numerators and denominators have obvious counting interpretations. We were unable to find such a proof and take a different route.

\subsection{Proof of Lemma~\ref{lem:betbeatsalp}}  That $f(m,k)\ge 1$ for all $1\le k\le m$, follows from three assertions:
\benu
\item $f(m,k+1)\ge f(m,k)$ for any $1\le k\le m-1$ and any $m\ge 1$.
\item $f(m,1)\ge f(m+1,1)$ for all $m\ge 1$.
\item $f(m,1)\to 1$ as $m\to \infty$.
\eenu
From the second and third assertions, it follows that $f(m,1)\ge 1$ for all $m$ and then the first assertion gives $f(m,k)\ge f(m,1)$, completing the proof. Now we prove the three assertions.

\para{Step-1} For $1\le k\le m-1$, canceling many factorials, $f(m,k+1)/f(m,k)$ is seen to be equal to
\bas
 \frac{(m+k+1)(4m-2k)(4m-2k+1)}{(m-k)(4m+2k+1)(4m+2k+2)} 
&= \frac{(4m+4k+4)(4m-2k)(4m-2k+1)}{(4m-4k)(4m+2k+1)(4m+2k+2)}.
\eas
Write $x=4m$ and $a=2k$. If we subtract $(x-2a)(x+a+1)(x+a+2)$ (the denominator) from  $(x+2a+4)(x-a)(x-a+1)$ (the numerator), we get $2x(x-2a)+2x+4a^2(a+2)$. As $x-2a=4(m-k)>0$, all three summands are positive if $1\le k\le m-1$ and $m\ge 1$. Thus $f(m,k+1)\ge f(m,k)$.

\para{Step-2} Observe that $f(m+1,1)/f(m,1)$ is equal to 
\bas
&\frac{m(m+2)\times (8m+2)\ldots (8m+9)}{2^6\times (2m+1)(2m+2)\times (4m+3)\ldots (4m+6)\times (4m)\ldots (4m+3)} \\
&= \frac{(m+2)(8m+3)(8m+5)(8m+7)(8m+9)}{2^4(2m+1)(2m+2)(4m+3)(4m+5)(4m+6)}
\eas
by canceling five factors in the numerator and denominator. Multiplying by powers of $2$ to write each terms as $8m+[\cdot]$, we see that
\bas
\frac{f(m+1,1)}{f(m,1)}&= \frac{(8m+16)(8m+3)(8m+5)(8m+7)(8m+9)}{(8m+4)(8m+8)(8m+6)(8m+10)(8m+12)} \\
&= \frac{(x+16)(x+3)(x+5)(x+7)(x+9)}{(x+4)(x+6)(x+8)(x+10)(x+12)}
\eas
with $x=8m$. Then the denominator minus the numerator is (we used Mathematica here) $15 (528 + 257 x + 40 x^2 + 2 x^3)$ which is positive for $x>0$. Thus  $f(m+1,1)\le f(m,1)$ for all $m\ge 1$.
\berk\label{rem:schurconcavity} Here is a less computational way to check the last point. For $x>0$, the function $\psi(a_{1},\ldots ,a_{k}):= (x+a_{1})\ldots (x+a_{k})$ is Schur-concave. This can be seen either from the well-known fact that elementary symmetric polynomials are Schur-concave (and $\psi$ is a positive linear combination of those) or by directly checking the condition for Schur-concavity: $(a_{i}-a_{j})(\frac{\partial}{\partial a_{i}}-\frac{\partial}{\partial a_{j}})\psi\le 0$. In our case,  $(16,9,7,5,3)$ majorizes $(12,10,8,6,4)$, hence $\psi(16,9,7,5,3)\le \psi(12,10,8,6,4)$. 
\eerk

\para{Step-3} That $f(m,1)\to 1$ as $m\to \infty$ is clear from the local central limit theorem, but can also argue directly from Stirlings' approximation:
\bas
f(m,1) &=  \frac{ (8m+1)! \times (m+1)!\times (m-1)!}{2^{6m}\times (2m)!\times (4m+2)!\times (4m-1)! } \\
&=   \frac{ (8m)! \times m! \times m!}{2^{6m}\times (2m)!\times (4m)!\times (4m)!} \times \frac{(8m+1)(m+1)(4m)}{m(4m+1)(4m+2)} \\
&\sim \frac{(8m)^{8m+\frac12} m^{2m+1}}{2^{6m}(2m)^{2m+\frac12}(4m)^{8m+1}}\times 2 \\
&= 1.
\eas

\subsection{Proof of Lemma~\ref{lem:alpbeatsbet}}  
The proof of the first statement in the Lemma will be  achieved in three steps, similarly to the proof of Lemma~\ref{lem:betbeatsalp}. The second statement follows from the first by setting $h=1$.
\benu
\item $g_{h}(m,k+1)\le g_{h}(m,k)$  if $h+1\le k\le [(4h-1)m^{2}]^{\frac13}$.
\item $g_{h}(m+1,h+1)\ge g_{h}(m,h+1)$ for $m\ge h+1$.
\item $g_{h}(m,h+1)\to 1$ as $m\to \infty$, for fixed $h$.
\eenu
The second and third asertions show that $g_{h}(m,h+1)\le 1$ for all $m\ge h+1$ and the first assertion shows that $g_{h}(m,k)\le 1$ if  $h+1\le k\le [(4h-1)m^{2}]^{\frac13}$. The proofs of the assertions are similar to that of Lemma~\ref{lem:betbeatsalp}.

\para{Step-1} Canceling many terms, $g_{h}(m,k+1)/g_{h}(m,k)$ is seen to be equal to
\bas
& \frac{(m+k+1-h)(4m-2k)(4m-2k+1)}{(m-k+h)(4m+2k+1)(4m+2k+2)}  \\
&= \frac{(x+2a-4(h-1))(x-a)(x-a+1)}{(x-2a+4h)(x+a+1)(x+a+2)}
\eas
with $x=4m$ and $a=2k$. Subtracting the numerator from the denominator gives
\bas
2((4h-1)x^{2}-2a^{3})+2x(8h+2a-1)+8(a^{2}(h-1)+ah+h)
\eas
which is positive provided $2a^{3}\le (4h-1)x^{2}$. Thus,  $g_{h}(m,k+1)\le g_{h}(m,k)$ if $k^{3}\le (4h-1)m^{2}$.

\para{Step-2} Let $(x)_{\uparrow k}=x(x+1)\ldots (x+k-1)$ and $(x)_{\uparrow\!\! \uparrow k}=x(x+2)\ldots (x+2k-2)$. Consider 
\small
\ba
\frac{g_{h}(m+1,k)}{g_{h}(m,k)} &= \frac{(8m+2)_{\uparrow 8}\times (m+1+k-h)\times (m+1-k+h)}{2^{6}(2m+1)(2m+2)\times (4m+2k+1)_{\uparrow 4}\times (4m-2k+2)_{\uparrow 4}} \\
&= \frac{(8m+2)_{\uparrow 8}\times (8m+8+8k-8h)\times (8m+8-8k+8h)}{(8m+4)(8m+8)\times (8m+4k+2)_{\uparrow\!\! \uparrow 4}\times (8m-4k+4)_{\uparrow\!\! \uparrow 4}}.
\ea
\normalsize
 Cancel $(8m+4)(8m+8)$ in the numerator and denominator, and set $x=8m$. The  above expression becomes $\prod_{j=1}^{8}(x+t_{j})/\prod_{j=1}^{8}(x+s_{j})$ where
\small
\ba
t &=(8k-8h+8,9,7,6,5,3,2,-8k+8h+8), \\
s &=(4k+8,4k+6,4k+4,4k+2,-4k+10,-4k+8,-4k+6,-4k+4).
\ea
\normalsize
As $k\ge h+1\ge 2$, these vectors are written in decreasing order. The vector of partial sums of $s-t$ is 
\small
\ba
 (8h-4k,8h-3,8h+4k-6,8h+8k-10,8h+4k-5,8h,8h-4k+4,0)
\ea
\normalsize
which is non-negative if and only if $k\le 2h$. Thus, for $h+1\le k\le 2h$, we see that $s$ majorizes $t$, and  (see Remark~\ref{rem:schurconcavity}) by Schur concavity $\prod_{j}(x+t_{j})\ge \prod_{j}(x+s_{j})$. That is, $g_{h}(m+1,k)\ge g_{h}(m,k)$.  In particular, this always holds for $k=h+1$.

\para{Step-3} $g_{h}(m,k)\to 1$ as $m\to \infty$, for any fixed $k,h$. This is immediate from Stirlings's formula or the local central limit theorem. Alternately one may observe that
\bas
g_{h}(m,k)= f(m,k)\frac{(m-k+1)\ldots (m-k+h)}{(m+k-h+1)\ldots (m+k)}
\eas
and use that $f(m,k)\to 1$ as $m\to \infty$.

\section{The comparison between the tails of $2S_{2m}$ and $S_{8m}$}\label{sec:s2mands8mtails}
We now prove the following crucial lemma comparing the tails of the two random variables.
\begin{lemma}\label{lem:s2mands8mtails} For some $m_{0}$ and any $m\ge m_{0}$, we have
\benu
\item $\overline{\alp}_{m}(k)\le \overline{\bet}_{m}(k)$ for all $k\ge 1$.
\item $\overline{\alp}_{m}(k)\ge \overline{\bet}_{m}(\ell)$ if $1\le k\le \ell-\frac{\ell^{2}}{4m}-1$.
\eenu
\end{lemma}
In Lemma~\ref{lem:alpbeatsbet}, we had the relationship $k=\ell-C\frac{\ell^{3}}{m^{2}}$. If the same relationship could be carried over to Lemma~ \ref{lem:s2mands8mtails}, we would have ended up with the Carter-Pollard~\cite{carterpollard} improvement of Tusn\'{a}dy's lemma. However, the proof (of the second part) got way longer, and hence we settled for the weaker form with $\ell^{2}/m$.

In addition to Lemma~\ref{lem:alpbeatsbet} and Lemma~\ref{lem:betbeatsalp}, we need the following basic estimates for binomial coefficients. Introduce the notation $D(p):=p\log(2p)+(1-p)\log(2-2p)$, usually written as  $D(\mb{Ber}(p)\|\mb{Ber}(1/2))$, for the relative entropy of $\mb{Bernoulli}(p)$ with respect to $\mb{Bernoulli}(1/2)$.
\benu
\item Lemma~4.7.1 of Ash~\cite{ash} states that with $A=1/\sqrt{8}$ and $B=1/\sqrt{2\pi}$, for any $1\le k\le n-1$,
\small
\begin{align}\label{eq:binestptwise}
A\frac{\sqrt{n}}{\sqrt{k(n-k)}}e^{-nD(\frac{k}{n})} \; \le \; \frac{1}{2^{n}}\binom{n}{k} \; \le \; B\frac{\sqrt{n}}{\sqrt{k(n-k)}}e^{-nD(\frac{k}{n})}.
\end{align}
\normalsize
\item Lemma~4.7.2 of Ash~\cite{ash} states that for $\frac{n}{2}<k<n$, with $A=1/\sqrt{8}$, 
\small
\begin{align}\label{eq:binesttail}
A\frac{\sqrt{n}}{\sqrt{k(n-k)}}e^{-nD(\frac{k}{n})} \; \le \; \frac{1}{2^{n}}\sum_{j=k}^{n}\binom{n}{j} \; \le  \; e^{-nD(\frac{k}{n})}.
\end{align}
\normalsize
\eenu
Both these estimates are proved using Stirlings formula, but one needs more than the first term in the asymptotic expansion. For a somewhat simpler proof (that uses only the second  term in Stirling's formula) of \eqref{eq:binestptwise} with the weaker constant $B=1/\sqrt{\pi}$, see Lemma~17.5.1 of  Cover and Thomas~\cite{coverthomas}. 

\para{The relative entropy function}  Let $s=2t-t^{2}$, so that $t\mapsto s$ is an increasing bijection of $[0,1]$ with itself. Define $Q(t):=4D\l(\frac12+\frac{t}{2}\r) -D\l(\frac12+\frac{s}{2}\r)$, an object that will occur repeatedly in the proof. Note that  $Q(t)=2t^{3}+O(t^{4})$ as $t\to 0$ and that $Q(t)$ is strictly positive for all $t>0$. Therefore, for any $\eps>0$, we can write $Q(t)\ge (2-\eps)t^{3}$ for $0\le t\le \del$ for some $\del>0$ and $Q(t)\ge c$ for $\del\le t\le 1$ for some $c>0$. For the sake of simplifying the writing, we use the more convenient form
\small
\begin{align}\label{eq:cubicboundforQ}
Q(t) \ge \frac32 t^{3} \;\; \mb{ for } \ 0\le t\le 1,
\end{align}
\normalsize
which is easiest to  check on a software, numerically or symbolically (in fact $3/2$ can be replaced by $1.65\cdots$). 


\berk The precise constants do not matter in the big picture. If $A$ were smaller or $B$ were larger, or if the factor in front of $m^{2/3}$ in Lemma~\ref{lem:alpbeatsbet} were less than $3^{1/3}$, the proof would still go through,  but with a choice of $s=t-\mu t^{2}$ for a sufficiently large $\mu$ and appropriate modifications that follow from it. This would only affect the constants in Lemma~\ref{lem:snandzcoupling} and Theorem~\ref{thm:snands4ncoupling}. The clean choice $\mu=1$ happens to work, hence we fix it up front. See Section~\ref{sec:bestconstant} for the best constant this proof can give.
\eerk

\medskip
With these preparations, we begin the proof of Lemma~\ref{lem:s2mands8mtails}. Again we assume that $n=2m$ is even. For brevity of notation,  write $\alp_{k},\bet_{k}$ for $\alp_{m}(k),\bet_{m}(k)$ and define the tails, $\overline{\alp}_{k}=\sum_{j\ge k}\alp_{j}$ and $\overline{\bet}_{k}=\sum_{j\ge k}\bet_{j}$. Recall that $\alp_{k}=\P\{S_{2m}=2k\}$ and $\bet_{k}=\P\{S_{8m}=4k-2\mb{ or }4k\}$  and hence
\small
\begin{align}\label{eq:tailsofalpbet}
\overline{\alp}_{k}=\frac{1}{2^{2m}}\sum_{j\ge k}\binom{2m}{m+j} \qquad \mb{ and }\qquad \overline{\bet}_{k}=\frac{1}{2^{8m}}\sum_{j\ge 2k-1}\binom{8m}{4m+j}.
\end{align}
\normalsize

\bprf[Proof of the first part of Lemma~\ref{lem:s2mands8mtails}] 
Lemma~\ref{lem:betbeatsalp} immediately implies that $\overline{\alp}_{k}\le \overline{\bet}_{k}$ if $k\ge 1$.
\eprf

%

\bprf[Proof of the second part of Lemma~\ref{lem:s2mands8mtails}] We prove it in stages, starting with the larger values of $\ell$ and proceeding to smaller values. We shall assume that $m$ is sufficiently large (so that expressions like $Cm^{2/3}\le m$ hold) without further comment.


\noindent{\bf Case $2(m^{2}\log m)^{1/3}\le \ell< 2m$:} 
Let $t=\frac{\ell}{2m}$ and $s=2t-t^{2}$. If $k\le \ell-\frac{\ell^{2}}{4m}$, then $\frac{k}{m}\le s$. Therefore, by \eqref{eq:tailsofalpbet} and \eqref{eq:binesttail} 
\bas
\frac{\overline{\alp}_{k}}{\overline{\bet}_{\ell+1}} &\ge \frac{A\sqrt{2m}}{\sqrt{(m+k)(m-k)}}e^{2m[4D(\frac12+\frac{\ell}{4m})-D(\frac12+\frac{k}{2m})]} \\
&\ge \frac{1}{2\sqrt{m}}e^{2mQ(t)} 
\eas
as $A=1/\sqrt{8}$ and $(m-k)(m+k)\le m^{2}$. By the bound~\eqref{eq:cubicboundforQ} $Q(t)\ge \frac32 t^{3}$, we see that   $2mQ(t)\ge 2\log m$ and hence the above expression is more than $\frac12 m^{3/2}$. Thus $\overline{\alp}_{k}\ge \overline{\bet}_{\ell+1}$.

\smallskip
\noindent{\bf Case $(\frac83 m^{2})^{1/3}\le \ell \le 2(m^{2}\log m)^{1/3}$:}
Again let $t=\frac{\ell}{2m}$ and $s=2t-t^{2}$ and observe that the condition $k\le \ell-\frac{\ell^{2}}{4m}$ ensures that $\frac{k}{m}\le s$. Then by  \eqref{eq:binestptwise}, 
\bas
\frac{\alp_{k}}{\bet_{\ell+1}} &\ge  \frac{A}{B}\frac{4m-2\ell+1}{8m+1}\frac{\sqrt{(4m+2\ell)(4m-2\ell)}}{2\sqrt{(m+k)(m-k)}}e^{2m[4D(\frac12+\frac{\ell}{4m})-D(\frac12 +\frac{k}{2m})]} \\
&\ge \frac{\sqrt{3}A}{4B}e^{2mQ(t)} 
\eas
if $m$ is large enough that $\ell\le m$, since in that case $(4m-2\ell)(4m+2\ell)\ge 12m^{2}$ and $(m-k)(m+k)\le m^{2}$ and $\frac{4m-2\ell+1}{8m+1}\ge \frac14$ when $\ell\le m$ (which holds for $\ell$ in this range as we assume $m$ is large enough). Plugging in the values of $A,B$ and using the bound \eqref{eq:cubicboundforQ}, we see that for all $\ell$ in this range
\small
\begin{align}\label{eq:ptbdforalpoverbet}
\alp_{k}\ge \frac{\sqrt{3\pi}}{8} e^{\frac{3\ell^{3}}{8m^{2}}}\bet_{\ell+1} \; \mb {whenever }1\le k\le \ell-\frac{\ell^{2}}{4m}.
\end{align}
\normalsize

Now suppose $1\le \ell \le L:=\lceil 4(m^{2}(\log m))^{1/3} \rceil$. Since we have shown that $\overline{\alp}(L+1)\ge \overline{\bet}(L+1)$ in the previous case,  it suffices to prove that $\sum_{j=k}^{K}\alp_{j}\ge \sum_{j=\ell}^{L}\bet_{j}$ whenever $k\le \ell-\mu\frac{\ell^{2}}{m}-1$ and $K=\lfloor L- \frac{L^{2}}{4m}\rfloor $. 

Let $\phi:\R_{+}\mapsto \R_{+}$ be a continuous decreasing function such that $\phi(k)=\alp_{k}$ for integer $k\ge 0$. Let $\kappa(\ell)=\ell-\frac{\ell^{2}}{4m}$ (not necessarily an integer). From \eqref{eq:ptbdforalpoverbet} we also have that 
\bas
\sum_{j=\ell}^{L}\phi(\kappa(j)) \; \ge \; \frac{\sqrt{3\pi}}{8} e^{\frac{3\ell^{3}}{8m^{2}}} \sum_{j=\ell}^{L}\bet_{j}.
\eas
Although $\kappa(j)$ ranges between $k$ and $K$, the sum on the left side cannot be simply bounded by $\sum_{i=k}^{K}\alp_{i}$, since a particular index $i$ may occur  as $\lfloor \kappa(j)\rfloor$ for more than one $j$ between $\ell$ and $L$. For $j$ in this range, $\kappa(j+1)-\kappa(j)\ge \eta:=1-\frac{2L+1}{m}$. As $\phi$ is decreasing, 
\bas
\sum_{j=\ell+1}^{L}\phi(\kappa(j)) &\le \frac{1}{\eta} \sum_{j=\ell+1}^{L}\phi(\kappa(j))\times (\kappa(j)-\kappa(j-1)) \\
&\le  \frac{1}{\eta} \int_{\kappa(\ell)}^{\kappa(L)}\phi(x)dx \\
&\le  \frac{1}{\eta}\sum_{i=k}^{K}\phi(i).
\eas
In the last line we bounded the integral by the sum, valid since $k\le \kappa(\ell)$ and $K\ge \kappa(L)$ are integers (and the spacings are of unit length). Plugging in the expression for $\eta$ and using $\phi(i)=\alp_{i}$,
\bas
\sum_{i=k}^{K}\alp_{i} &\ge \frac{\sqrt{3\pi}}{8}\l(1-\frac{2L+1}{m}\r)e^{\frac{3\ell^{3}}{8m^{2}}}\sum_{j=\ell}^{L}\bet_{j}. 
\eas
For $m$ large, the product of the first two factors is more than $1/e$, hence it suffices to have $\ell^{3}\ge \frac{8}{3}m^{2}$ to conclude that $\overline{\alp}_{k}\ge \overline{\bet}_{\ell}$.

\smallskip
\noindent{\bf Case $1\le \ell\le (\frac83 m^{2})^{1/3}$:} Let $L':=(\frac83 m^{2})^{1/3}$ and $L''=(3m^{2})^{1/3}$ and fix $\ell\le L'$. From the second statement in Lemma~\ref{lem:alpbeatsbet}, we know that $\bet_{\ell+2}+\bet_{\ell+2}+\ldots +\bet_{L''}$ can be bounded above by $\alp_{\ell+1}+\alp_{\ell+2}+\ldots +\alp_{L''-1}$.   We claim that the entire tail $\overline{\bet}_{L''+1}$ can be bounded by one term $\alp_{\ell}$. When added to the above inequality, this  proves that $\overline{\alp}_{\ell}\ge \overline{\bet}_{\ell+2}$. This completes the proof, since $k\le \ell-2$.

To prove the claim, observe that by \eqref{eq:binesttail} and \eqref{eq:binestptwise}, and the bounds for $D$, 
\bas
\overline{\bet}_{L''+1}&\le e^{-8mD(\frac12+\frac{L''}{4m})} \; \le \; e^{-\frac{(L'')^{2}}{m}}, \\
\alp_{\ell}&\ge \frac{A\sqrt{2m}}{\sqrt{(m-k)(m+k)}}e^{-2mD(\frac12+\frac{\ell}{2m})}  \; \ge  \; \frac{1}{2\sqrt{m}}e^{-\frac{\ell^{2}}{m}-2E\frac{\ell^{4}}{m^{3}}}.
\eas
As $\ell\le L'$ and $L''$ is larger than $L'$ by a factor more than $1$, it is clear that $\alp_{\ell}\ge \overline{\bet}_{L''+1}$ for $\ell \le L'$.
\eprf

\section{Proof of Theorem~\ref{thm:snands4ncoupling}}\label{sec:binomialcouplingandgeneralcoupling} 
Finally we deduce Theorem~\ref{thm:snands4ncoupling} from Lemma~\ref{lem:s2mands8mtails}.  The basic idea is simple and stated as follows. 
\begin{lemma}\label{lem:gencoupling} Let $\alp=(\alp_{k})_{k\ge 0}$ and $\bet=(\bet_{k})_{k\ge 0}$ be probability mass functions. Let $\overline{\alp}(x)=\sum_{j\ge x}\alp_{j}$ and the similarly defined $\overline{\bet}$ denote their tails. Suppose that $\overline{\alp}(k-f(k))\ge\overline{\bet}(k)$ and $\overline{\bet}(k-g(k))\ge \overline{\alp}(k)$ for some $f,g:\N\mapsto \N$ and for all $k\ge 1$.  Then, there is a coupling of $X\sim \alp$ and $Y\sim \bet$ such that $X\ge Y-f(Y)$ and $Y\ge X-g(X)$.
\end{lemma}
Assuming this, let us deduce Theorem~\ref{thm:snands4ncoupling}.
\bprf[Proof of Theorem~\ref{thm:snands4ncoupling}]  Observe that what appear in Lemma~\ref{lem:s2mands8mtails} are (we drop the subscript $m$ on $\overline{\alp},\overline{\bet}$ to simplify notation) are $\overline{\alp}(k)=\P\{2S_{2m}\ge 4k\}$ and $\overline{\bet}(k)=\P\{S_{8m}\ge 4k-2\}$. Therefore, by the conclusion of that lemma, if $m$ is large enough, we get
\benu
\item $\P\{2S_{2m}\ge j\}\le \P\{S_{2m}\ge j-2\}$. To see this, observe that the left hand side equals $\overline{\alp}(k+1)$ for $4k+1\le j\le 4k+4$, and by Lemma~\ref{lem:s2mands8mtails}, it does not exceed $\P\{S_{8m}\ge 4k+2\}$ which is at most $\P\{S_{8m}\ge j-2\}$ as $j-2\le 4k+2$.
\item $\P\{S_{8m}\ge j\}\le \P\{2S_{2m}\ge j-1-4\theta(j/4)\}$ where $\theta(x)=\frac{x^{2}}{4m}+2$. To see this, we use the bound $\P\{S_{8m}\ge j\}\le \overline{\bet}(k)$, for $4k-2\le j\le 4k+1$. By Lemma~\ref{lem:s2mands8mtails}, this is at most $\P\{2S_{2m}\ge 4(k-\theta(k))\}$. Then we observe that in each of the four cases, $4(k-\theta(k))$ is at least as large as  $j-1-4\theta(j/4)$.
\eenu
By the symmetry of $S_{2m}$ and $S_{8m}$, the above probability comparisons hold for the absolute values. Lemma~\ref{lem:gencoupling}, assures that $|S_{2m}|$ and $|S_{8m}|$ can be coupled so that
\bas
|S_{8m}|-9-\frac{|S_{8m}|^{2}}{16m} \le 2|S_{2m}|\le |S_{8m}|+2.
\eas
Conditional on the absolute values, choose the same sign for both to get a coupling of $S_{2m}$ and $S_{8m}$ as in the statement of the theorem.
\eprf
Here is the obvious of the  lemma on coupling.
\bprf[Proof of Lemma~\ref{lem:gencoupling}] Let $V\sim \mb{uniform}[0,1]$, find $k,\ell$ such that $\overline{\alp}(k)\ge V>\overline{\alp}(k+1)$ and $\overline{\bet}(\ell)\ge V>\overline{\bet}(\ell+1)$ and set $X=k$ and $Y=\ell$. Then $X\sim \alp$ and $Y\sim \bet$. Further, $\overline{\alp}(k)>\overline{\bet}(\ell+1)$ and $\overline{\bet}(\ell)>\overline{\alp}(k+1)$. But by the assumptions, $\overline{\bet}(k-g(k))\ge \overline{\alp}(k)$ and $\overline{\alp}(\ell-f(\ell))\ge \overline{\bet}(\ell)$. Therefore, we see that $\overline{\bet}(k-g(k))>\overline{\bet}(\ell+1)$ and $\overline{\alp}(\ell-f(\ell))>\overline{\alp}(k+1)$ which of course imply that $k-g(k)\le \ell$ and $\ell-f(\ell)\le k$.
\eprf

\section{A remark on the optimal constant}\label{sec:bestconstant} 
What is the best constant that can be achieved by this proof? Fix $\mu\le 1$ and make the choice $s=(2t-\mu t^{2})\wedge 1$ (i.e., for $t\ge (1-\sqrt{1-\mu})/\mu$ we take $s=1$), and use the expansion $Q(t)=2\mu t^{3}-O(t^{4})$. This is how various components of the proof change. Let $a_{m}=a(m^{2}\log m)^{1/3}$ and $b_{m}=(bm^{2})^{1/3}$ mark the boundaries  between the three cases in the proof of Lemma~\ref{lem:s2mands8mtails}. Below $t=\ell/2m$. The first case offers no problems, we simply use $Q(t)\ge c a_{m}^{3}$ for all $t\ge a_{m}$, and choose $a$ large enough  to ensure that $2mQ(a_{m})\ge 2\log m$. In the second case, we need $\frac{\sqrt{3\pi}}{8}e^{2mQ(t)}\ge 1$ for $b_{m}\le t\le a_{m}$, which only  needs to be checked at $\ell=b_{m}$. As $Q(t)\sim 2\mu t^{3}$, this translates to $b> \frac{2}{\mu}\log(8/\sqrt{3\pi})$. For this to meet up with the third case and cover the entire range of $\ell$, we require $b< 3$ (since the third case, that depends on Lemma~\ref{lem:alpbeatsbet}, goes up to $\ell \le (3m^{2})^{1/3}$). This gives $\mu> \frac23 \log(8/\sqrt{3\pi})=0.6385...$.

What does this mean for the constant in Lemma~\ref{lem:snandzcoupling}? With $\mu=1$, we got  $1\times Z^{2}$ in the statement and with $\mu$, we could get  $\mu Z^{2}$, following exactly the same proof. Actually we were wasteful in the proof of Lemma~\ref{lem:snandzcoupling} when we bounded $(Z+\frac{c}{\sqrt{n}})^{2}$ by $2Z^{2}+2\frac{c^{2}}{n}$ in \eqref{eq:wastefulboundonsquareofsum}, thereby losing a factor of $2$. We could instead use $(a+b)^{2}\le a^{2}/p+b^{2}/q$ for any positive $p,q$ that sum to $1$, with a $p$ close to $1$. This would  improve the bound in Lemma~\ref{lem:snandzcoupling} to $cZ^{2}+C$ for any $c>1/2$ with the choice $\mu=1$, and to $0.319...$ with the optimal choice of $\mu$. This still falls short of Tusnady's constant $1/4$, not to mention our suboptimal additive constants.

\newpage
\part{Coupling via Markov chains or Stein coefficient}\label{partII}

This part was inspired entirely by Chatterjee's proof~\cite{chatterjee}, in an attempt to refashion it to  compare Binomials among themselves than with the Gaussian. If a pair of probability distributions are stationary distributions of nearest neighbour Markov chains on integers, then we show how to construct a coupled Markov chain on $\Z^{2}$ whose stationary distribution provides a coupling of the two given distributions. We prove a general result that the coupling is good if the transition rates are close to each other (we write the result for the special class of {\em Ehrenfest-like chains} that are determined by one function on integers called a {\em Stein coefficient}). This is applied to suitable choices of  Markov chains to prove Theorem~\ref{thm:bingausscouplingcauchy} and Theorem~\ref{thm:hypergeogausscouplingcauchy}. 

\section{Coupling two nearest-neighbour Markov chains on integers}
Let $\mathbb S=\{a,a+1,\ldots ,b\}$ and $\mathbb T=\{c,c+1,\ldots ,d\}$ be shifted finite segments  of integers ($b-a$ and $d-c$ must be integers, but $a,b,c,d$ need not be). Let $X$ be a continuous-time, nearest-neighbour Markov chain on $\mathbb S$ with transitions $i\to i\pm 1$ with rates $\lam_{i}^{\pm}$. Similarly, let $Y$ be a nearest-neighbour chain on $\mathbb T$ with rates $\mu_{j}^{\pm}$. Naturally, $\lam_{b}^{+},\mu_{d}^{+},\lam_{a}^{-},\mu_{c}^{-}$ are all zero. We assume that the chains are irreducible and denote the unique stationary probability distributions of $X$ and $Y$ by  $\alp$ and $\bet$ respectively. 

\subsection{A coupled pair of Markov chains} The goal is to construct a Markov chain on $\mathbb U:=\mathbb S\times \mathbb T$ such that the co-ordinates move like $X$ and $Y$, but stay close to each other.  The most natural idea would be to make $X$ and $Y$ take a step to the right together or a step to the left together, to the extent possible, and  avoid taking steps in opposing directions. This leads us to the Markov chain $Z$ on $\mathbb U$ with rates  
\bas
\theta_{i,j}^{+,+}=\lam_{i}^{+}\wedge \mu_{j}^{+}, \qquad \theta_{i,j}^{+,\circ}=(\lam_{i}^{+}-\mu_{j}^{+})_{+}, \qquad \theta_{i,j}^{\circ,+}=(\mu_{j}^{+}-\lam_{i}^{+})_{+}, \\
\theta_{i,j}^{-,-}=\lam_{i}^{-}\wedge \mu_{j}^{-}, \qquad \theta_{i,j}^{-,\circ}=(\lam_{i}^{-}-\mu_{j}^{-})_{+}, \qquad \theta_{i,j}^{\circ,-}=(\mu_{j}^{-}-\lam_{i}^{-})_{+}.
\eas
The notation is self-explanatory: For example, $\theta_{i,j}^{+,\circ}$ is the transition rate from $(i,j)$ to $(i+1,j)$ while $\theta_{i,j}^{-,-}$ is the transition rate from $(i,j)$ to $(i-1,j-1)$. Observe that there are no transitions to $(i+1,j-1)$ or $(i-1,j+1)$. The generator $L$ of $Z$ acts on $f:\mathbb U\mapsto \R$ as
\bas
Lf(i,j) &= \theta_{i,j}^{+,+}[f(i+1,j+1)-f(i,j)] + \theta_{i,j}^{-,-}[f(i-1,j-1)-f(i,j)] \\
& \;\; +\theta_{i,j}^{+,\circ}[f(i+1,j)-f(i,j)] +\theta_{i,j}^{-,\circ}[f(i-1,j)-f(i,j)] \\
&\;\; +\theta_{i,j}^{\circ,+}[f(i,j+1)-f(i,j)] +\theta_{i,j}^{\circ,-}[f(i,j-1)-f(i,j)].
\eas
\subsection{The stationary distribution of $Z$ is a coupling of $\alp$ and $\bet$}
As $\mathbb U$ is finite,  $Z$ necessarily has  a stationary distribution $\gam$. Then $\E_{\gam}[Lf(Z)]=0$ for all $f:\mathbb U\mapsto \R$.  In fact we shall show shortly that there is a unique stationary distribution, but the uniqueness will not play a role in the  analysis.

Let $f(i,j)=\phi(i)$ and $g(i,j)=\psi(j)$. From the formula for $L$, we see that 
$Lf(i,j)=L_{1}\phi(i)$ and $Lg(i,j)=L_{2}\psi(j)$ where $L_{1}$ and $L_{2}$ are the generators of $X$ and $Y$, respectively. The equations $\E_{\gam}[Lf(Z)]=0$ and $\E_{\gam}[Lg(Z)]=0$ show that the co-ordinates of $Z$ move like $X$ and $Y$ and that the marginals of $\gam$ are $\alp$ and $\bet$. In other words, $\gam$ is a coupling of $\alp$ and $\bet$.

\para{Uniqueness of the stationary distribution} 
In general, and even in examples of interest to us, $Z$ is not irreducible. But we claim that it has a unique recurrent class, and hence a unique stationary probability distribution. To prove this claim, choose $M'$ and $M''$ to be medians of $\mathbb S$ and $\mathbb T$ in such a way that the line $\ell_{*}$ in $\R^{2}$ having slope $1$ and passing through $(M',M'')$ either intersects the top and bottom sides of $[a,b]\times [c,d]$ or the left and right sides of $[a,b]\times [c,d]$. This is trivial when the medians are unique, for then $(M',M'')$ is the center of the rectangle $[a,b]\times[c,d]$. A little case analysis shows that more generally, any choice of medians works, except when $\mathbb S$ and $\mathbb T$ have the same even number of elements, in which case two of the four choices (both being the smaller medians or both being the larger medians) work.

We claim that every state in $\mathbb U$ leads to $(M',M'')$. Why so? Within each diagonal line $\ell_{d}:=\{(i,j)\in \mathbb S\times \mathbb T\suchthat i-j=d\}$, every state leads to every other state.   Further, on the top edge (respectively bottom, right, left) of $\mathbb U$, there is a strictly positive rate to move to the right (respectively left, up, down). Therefore, if $(i,j)$ is to the ``left of'' $\ell_{*}$, then depending on whether $\ell_{*}$ intersects the bottom and top or left and right, move along the diagonal till you hit the  top or left sides, and then move right or down to get to $\ell_{*}$, and then move along $\ell_{*}$ to get to $(M',M'')$. A similar argument works if $(i,j)$ is to the right of the diagonal $\ell_{*}$. 

\subsection{Tail bounds on the difference $X-Y$ in the coupling $\gamma$}  Suppose $f(i,j)=\phi(i-j)$. Then writing $k=i-j$, we get
\bas
Lf(i,j) &= (\theta_{i,j}^{+,\circ}+\theta_{i,j}^{\circ,-})[\phi(k+1)-\phi(k)] + (\theta_{i,j}^{\circ,+}+\theta_{i,j}^{-,\circ})[\phi(k-1)-\phi(k)] \\
=&A(i,j)[\phi(k+1)-\phi(k-1)]+B(i,j)[\phi(k+1)-2\phi(k)+\phi(k-1)]
\eas
where 
\bas
A(i,j) &= \frac{\theta_{i,j}^{+,\circ}+\theta_{i,j}^{\circ,-}-\theta_{i,j}^{\circ,+}-\theta_{i,j}^{-,\circ}}{2}\; = \; \frac12 \l( \lam_{i}^{+}+\mu_{j}^{-}-\mu_{j}^{+}-\lam_{i}^{-}\r), \\
B(i,j) &= \frac{\theta_{i,j}^{+,\circ}+\theta_{i,j}^{\circ,-}+\theta_{i,j}^{\circ,+}+\theta_{i,j}^{-,\circ}}{2} \; = \; \frac12\l(|\lam_{i}^{+}-\mu_{j}^{+}|+|\lam_{i}^{-}-\mu_{j}^{-}|\r).
\eas
Let $Z=(X,Y)$ and $H=X-Y$. The equation $\E_{\gam}[Lf(Z)]=0$ becomes
\small
\begin{align*}
\E_{\gam}[A(Z)(\phi(H-1)-\phi(H+1))] = \E_{\gam}[B(Z)(\phi(H+1)-2\phi(H)+\phi(H-1))].
\end{align*}
\normalsize
Let $\psi(x)=\phi(x+1)-\phi(x)$ and rewrite this as 
\small
\begin{align}\label{eq:stateqforfunctionofdifference2}
\E_{\gam}[(B-|A|)\times (\psi(H)-\psi(H-1))] &= 2\E_{\gam}[A_{-}\psi(H-1)-A_{+}\psi(H)].
\end{align}
\normalsize
Here and below, $A_{+}=\max\{A,0\}$ and $A_{-}=(-A)_{+}$.

The basic idea now is to plug in various test functions to get bounds on the tail of the distribution of $H$ under $\gam$. That is precisely what we want, a coupling of $X$ and $Y$ so that $X-Y$ has light tails.  Although it is possible to do this in general and write some bounds, we specialize to the class of Markov chains that we use in this paper.

\subsection{Ehrenfest-like chains} Assume that $X$ and $Y$ have rates $\lam_{i}^{\pm}=S(i)\mp i$ and $\mu_{j}^{\pm}=T(j)\mp j$ for some $S:\mathbb S\mapsto \R_{+}$ and $T:\mathbb T\mapsto \R_{+}$.  When $\mathbb S=\{-\frac12 N,-\frac12 N+1,\ldots ,\frac12 N-1,\frac12 N\}$ and $\lam_{i}^{\pm}=\frac12N\mp i$, this gives the usual Ehrenfest chain, which is why we refer to these more general chains as Ehrenfest-like. For Ehrenfest-like chains,
\bas
A(i,j) &= j-i \\
B(i,j) &= \frac12 \l(|S(i)-T(j)+j-i|+|T(j)-S(i)+j-i|\r) \\
&= |j-i| + (|T(j)-S(i)|-|j-i|)_{+}
\eas
using $|x-y|+|x+y|=2|x|+2(|y|-|x|)_{+}$ (which is also $2\max\{|x|,|y|\}$). Thus, $A(Z)=-H$ and $B(Z)-|A(Z)|=(Q-|H|)_{+}$ where $Q=|T(Y)-S(X)|$. Consequently, \eqref{eq:stateqforfunctionofdifference2} becomes
\small
\begin{align}\label{eq:stateqforfunctionofdifferenceehrenfest}
\E_{\gam}\l[(Q-|H|)_{+}(\psi(H)-\psi(H-1))\r] = 2\E_{\gam}\l[H_{+}\psi(H-1)-H_{-}\psi(H)\r].
\end{align}
\normalsize
By plugging in various test functions, one gets bounds on expectations of functions of $H$ in terms of those of $Q$. To illustrate this, in \eqref{eq:bdfromMCcouplingforEhrenfest} and \eqref{eq:bdfromMCcouplingforEhrenfest2} we obtain a tail bound and a bound on expectations of functions in a straightforward way. But there is one place where it becomes necessary to have the tighter bound  \eqref{eq:bdfromMCcouplingforEhrenfest3}, hence we use that everywhere. 

\para{Tail bounds} Fix an integer $a\ge 0$ and let $\psi(x)=\one_{x\ge a}$. Then  \eqref{eq:stateqforfunctionofdifferenceehrenfest} leads to
$2\E_{\gam}[H_{+}\one_{H\ge a+1}] \le \E_{\gam}[(Q-a)_{+}]$. Similarly for $H_{-}$. Add the two to get
\small
\begin{align}\label{eq:bdfromMCcouplingforEhrenfest}
\P_{\gam}\{|H|\ge a+1\}  \; \le \;   \frac{1}{a+1} \E_{\gam}[|H|\one_{|H|\ge a+1}] \; \le \;  \frac{1}{a+1}\E_{\gam}[(Q-a)_{+}].
\end{align}
\normalsize

\para{Expectation bounds} Let $g:\Z \mapsto \R_{+}$ be increasing with $g(0)=0$. Let $\psi(x)=\frac{1}{x+1}g(x+1)\one_{x\ge 0}$. Then $\psi(x)-\psi(x-1)$ is at most $g(x+1)-g(x)$ (which is zero for negative $x$). Therefore \eqref{eq:stateqforfunctionofdifferenceehrenfest} implies that
\bas
2\E_{\gam}[g(H)] \; \le \; \E_{\gam}\l[(Q-|H|)_{+}\l(g(H+1)-g(H)\r) \r] \; \le \; \E_{\gam}[Qg(Q)].
\eas
Similarly for $-H$. Adding the two gives $\E_{\gam}[g(|H|)] \le \E_{\gam}[Qg(Q)]$. 

In particular, for   $g(x)=(e^{\theta x}-1)\one_{x\ge 0}$ this gives
\small
\begin{align}\label{eq:bdfromMCcouplingforEhrenfest2}
\E_{\gam}[e^{\theta|H|}] &\le 1+\E_{\gam}\l[Q(e^{\theta Q}-1)\r].
\end{align}
\normalsize
Chatterjee's method requires $\theta^{2}$ in the exponent on the right side, which motivates the more involved analysis that follows.

\para{Exponential moment bounds} Let $\psi(x)=e^{\theta (x+1)}\one_{x> 0}$ where $0\le \theta\le 1$. From \eqref{eq:stateqforfunctionofdifferenceehrenfest},
\bas
2\E_{\gam}[H_{+}e^{\theta H_{+}}] &= (e^{\theta}-1)\E_{\gam}[(Q-|H|)_{+}e^{\theta H}\one_{H> 0}].
\eas
Adding it to the corresponding inequality holds for $H_{-}$ gives
\small
\begin{align}\label{eq:prelimbdforHexpH}
2\E_{\gam}[|H|e^{\theta|H|}] &\le (e^{\theta}-1)\E_{\gam}[Qe^{\theta |H|}].
\end{align}
\normalsize
We  bounded $(Q-|H|)_{+}$ by $Q$.

We now use the fact that for any convex function $\phi$ and its Legendre or convex dual $\phi^{*}$ (both defined on subintervals of $\R$),  $xy\le \phi(x)+\phi^{*}(y)$, and hence $\E[XY]\le \E[\phi(X)]+\E[\phi^{*}(Y)]$ for any random variables $X,Y$ taking values in the domains of $\phi$ and $\phi^{*}$ respectively.  One may recall that this is indeed the idea behind H\"{o}lder's inequality, where one takes $\phi(x)=x^{p}/p$ and $\phi^{*}(x)=x^{q}/q$ for conjugate exponents $p,q$. But we apply it to a function $\phi$ smaller than $x^{p}$ for all $p>1$. 

Let  $\phi(x):=\bet x\log x$ for $x>0$, in which case  $\phi^{*}(y):=\bet e^{-1+\frac{y}{\bet}}$ for $y\in \R$. Applying the inequality to $X=e^{\theta |H|}$ and $Y=Q$ on the right side of \eqref{eq:prelimbdforHexpH}, 
\bas
2\E_{\gam}[|H|e^{\theta|H|}] &\le \bet\theta(e^{\theta}-1)\E_{\gam}[|H|e^{\theta |H|}]+\bet(e^{\theta}-1)e^{-1}\E_{\gam}[e^{Q/\bet}]).
\eas
Rearranging, we get
\bas
\E_{\gam}[\theta |H|e^{\theta|H|}] &\le \frac{\bet \theta(e^{\theta}-1)}{e(2-\bet\theta(e^{\theta}-1))} \E_{\gam}[e^{Q/\bet}].
\eas
as long as $\bet\theta(e^{\theta}-1)<2$. Writing  $\bet\theta(e^{\theta}-1)=2(1-\del)$ with $0<\del<1$,  (as $\theta\mapsto \theta(e^{\theta}-1)$ is a bijection from $\R_{+}$ to itself, for every $\theta>0$ and $\del\in (0,1)$, there is a unique $\bet>0$ such that $\bet\theta(e^{\theta}-1)=2(1-\del)$) we get
\bas
\E_{\gam}[\theta |H|e^{\theta|H|}] &\le \frac{1-\del}{e\del} \E_{\gam}\l[e^{\frac{1}{2(1-\del)}\theta(e^{\theta}-1)Q}\r] \\
&\le \frac{1-\del}{e\del} \E_{\gam}[e^{\frac{1}{2(1-\del)}e^{\theta}\theta^{2}Q}]
\eas
since $e^{\theta}-1=e^{\theta'}\theta$ for some $\theta'\in (0,\theta)$. Now for any $\mu>0$, we can write $e^{x}\le \frac{1}{\mu}xe^{x}+e^{\mu}$ for $x\ge 0$ (the first term suffices for $x\ge\mu$ and the second for $x\le \mu$). Therefore, 
\small
\begin{align}
\E_{\gam}[e^{\theta|H|}] &\le e^{\mu}+\frac{1-\del}{\mu \del e}\E_{\gam}\l[e^{\frac{1}{2(1-\del)}e^{\theta}\theta^{2}Q}\r] \nonumber \\
&\le \l(e^{\mu}+\frac{1-\del}{\mu \del e}\r)\E_{\gam}\l[e^{\frac{1}{2(1-\del)}e^{\theta}\theta^{2}Q}\r] \label{eq:bdfromMCcouplingforEhrenfest3}
\end{align}
\normalsize
valid for any $0<\del<1$ and $\mu>0$ and $\theta>0$. This is better than \eqref{eq:bdfromMCcouplingforEhrenfest2} when $\theta$ is small because of the power $2$ on $\theta$ in the exponent.  In examples, the right side is finite only for small $\theta$, hence \eqref{eq:bdfromMCcouplingforEhrenfest3} is actually the more useful one. 

In applying \eqref{eq:bdfromMCcouplingforEhrenfest3}, one can optimize over the right side, but it is perhaps more illuminating to observe the following.
\benu
\item Given $\eps>0$ and $\Theta<\infty$, choosing $\mu$ and $1-\del$ small, we see that there is some $M=M(\eps,\Theta)$ such that 
\small
\begin{align} \label{eq:bdfromMCcouplingforEhrenfest4}
\E_{\gam}[e^{\theta|H|}] &\le (1+\eps)\E_{\gam}\l[e^{M\theta^{2}Q}\r] \;\;\mb{ for }\theta\le \Theta.
\end{align}
\normalsize
\item Given $\eps>0$, choosing $\del$ small and  $\mu$ arbitrarily, we see that there is some $\theta_{\eps}>0$ and $M_{\eps}<\infty$ such that
\small
\begin{align} \label{eq:bdfromMCcouplingforEhrenfest5}
\E_{\gam}[e^{\theta|H|}] &\le M_{\eps}\E_{\gam}\l[e^{(\frac12+\eps)\theta^{2}Q}\r] \;\;\mb{ for }\theta\le \theta_{\eps}.
\end{align}
\normalsize 
\eenu

\berk Can one prove a Tusn\'{a}dy type lemma by this method? It is not possible in the naive sense, as the coupling $\gam$ is supported on a union of diagonals of $\mathbb U$. A more plausible approach is to interpret the inequalities between expectations that we have obtained as a stochastic domination result such as $|H|\prec aQ+b$ (roughly speaking). That would imply a coupling of $|H|$ and $Q$ so that $|H|\le aQ+b$. This can be done.

 However, it is not clear to us that one can couple variables $X,Y$ having the marginal $\alp,\bet$ so that $|X-Y|\le |T(Y)-S(X)|+1$. If this could be achieved in  some way, that would be nice because, in  the main example of coupling $2S_{n}$ and $S_{4n}$, the variable $Q=|T(Y)-S(X)|$ turns out to be exactly $O(1+\frac{1}{n}S_{4n}^{2})$, as required by Theorem~\ref{thm:snands4ncoupling}. 
\eerk

\section{Ehrenfest-like chains and the Stein coefficient}
To get couplings of given distributions on integers, we construct Ehrenfest-like Markov chains with these distributions and then use the general results of the previous section. The given distributions may need to be perturbed first to be able to do this. In this section we investigate Ehrenfest-like chains in general, although some of the discussions here are not needed in the proofs of Theorems~\ref{thm:bingausscouplingcauchy} and \ref{thm:hypergeogausscouplingcauchy}.
\subsection{Ehrenfest-like chain with given stationary distribution} Let $\alp$ be a probability distribution whose support $\mathbb S=\{a,a+1,\ldots ,b\}$ is a shifted segment of integers ($a,b$ need not be integers).   A nearest-neighbour Markov chain on $\mathbb S$ with transition rates $\lam_{i}^{\pm}$ is reversible for $\alp$ if and only if $\alp(i)\lam_{i}^{+}=\alp(i+1)\lam_{i+1}^{-}$ for all $i\in \mathbb S$ (with $\alp(i)$ and $\lam_{i}^{\pm}$ set to zero for $i\not\in \mathbb S$). Since we assume that $\alp(i)>0$ for all $i\in \mathbb S$, it is always possible to find such rates, for example, set  $\lam_{i}^{\pm}=\alp(i\pm 1)$.

For the chain to be Ehrenfest-like, the condition is that $\lam_{i}^{\pm}=T(i)\mp i$ for a function $T:\mathbb S\mapsto \R$. Does such a function exist? Unlike the general situation above where we had two parameters $\lam_{i}^{\pm}$ per state to play with, now we have only one parameter $T(i)$. The equations for reversibility take the form 
\small
\begin{align}\label{eq:Steinequationforbinomial}
\alp(i)(T(i)-i)=\alp(i+1)(T(i+1)+(i+1)).
\end{align}
\normalsize
This of course implies that $T(i)=|i|$ when $i\in \{a,b\}$, an end-point of $\mathbb S$. For $i\not\in \mathbb S$, since $\alp(i)=0$ the above equations say nothing about $T(i)$, and it need not be defined. Our convention is to set $T(i)=|i|$ for $i\not\in \mathbb S$.

Now assume that $\mathbb S=\{a,a+1,\ldots ,b\}$ is finite. Then we must set $T(b)=b$ and then successively solve for $T(b-1)$, $T(b-2)$, \dots ,$T(a)$. This is possible as $\alp(i)>0$ for all $i\in \mathbb S$. Inductively, the solution is seen to be 
\small
\begin{align}\label{eq:generalformulaforSteincoeff}
T(i)=i+\frac{2}{\alp(i)}\sum_{j: \ j>i}\alp(j)j \qquad \mb{ for }i\in \mathbb S.
\end{align}
\normalsize
But there is one more equation, namely \eqref{eq:Steinequationforbinomial} for $i=a-1$, that forces  $T(a)=-a$. This is  satisfied by the expression in \eqref{eq:generalformulaforSteincoeff} if and only if  $\alp$ has zero mean. In short, stationary distributions of Ehrenfest-like chains are precisely those that have zero mean and support equal to a shifted segment of integers. If $\alp$ is such a distribution, then there is a unique Ehrenfest-like chain that keeps it stationary. The function $T$ will be called the {\em Stein coefficient}\footnote{To reconcile with the language of Chatterjee's paper~\cite{chatterjee}, we could call it {\em Stein coefficient with respect to the Binomial distribution}, and what he uses as {\em Stein coefficient with respect to the Gaussian}.  The connection to Stein's method is this: $\E[\sig^{2}f'(W)]=\E[Wf(W)]$ for a rich class of functions $f$ if and only if $W\sim N(0,\sig^{2})$ and $\E[(n-X)f(X+1)-Xf(X)]=0$ for a rich class of functions if and only if $X\sim \mb{Binomial}(n,\frac12)$. Chatterjee's Stein coefficient for a random variable $W$ is a function $T$ that replaces $\sig^{2}$ in the first equation and our Stein function is a function $T$ that replaces $n$ in the second. The essence of Stein's method is that if $T$ is close to a constant, then the corresponding $W$ (respectively $X$) has a distribution close to Gaussian (respectively Binomial).} of $\alp$.

After presenting two examples that will be of use to us later, in the next two subsections we study the behaviour of Stein coefficient under convolution and scaling. The relevance is of course that we shall want to compare $S_{4n}$ with $2S_{n}$.

\beg \label{eg:binomialstein} Let $\alp$ be the distribution of the centered Binomial$(n,p)$ distribution. The support is $\mathbb S=\{-np,-np+1,\ldots ,nq\}$ where $q=1-p$, and the mass function is $\alp(x)=\binom{n}{j}p^{j}q^{n-j}$ where $j=x+np$.  The Stein coefficient of $\alp$ is $T(x)= 2pqn+(q-p)x$. Indeed,  the ratio of the right side to the left side of \eqref{eq:Steinequationforbinomial} is (again $x=j-np$)
\bas
\frac{2pqn+(q-p)(x+1)+(x+1)}{2pqn+(q-p)x-x}\times\frac{(n-j)p}{(j+1)q} = 1.
\eas
In particular, for the symmetric case $p=\frac12$ (which is the distribution of $S_{n}/2$),  we have $T(x)=\frac12 n$ for $x\in \{-\frac12 n,-\frac12n+1,\ldots ,\frac12 n-1, \frac12 n\}$. The corresponding chain is the usual Ehrenfest chain. Binomials are the only distributions for which the Stein coefficient is linear (this is the Stein characterizing equation for the Binomial).
\eeg
\beg \label{eg:hypergeometricstein} Let $\hat{S}_{k}[n,s]:=\frac12(S_{k}[n,s]+k)$ (see the text preceding Lemma~\ref{lem:hypergeogausscoupling} for the definition of $S_{k}[n,s]$). Equivalently, $\hat{S}_{k}$ is the sum of $k$ coupons drawn without replacement from a box containing $np=\frac12(n+s)$ coupons labeled $1$ and $nq=\frac12(n-s)$ coupons labeled $0$. Let $\alp$ denote the distribution of the centered variable $\hat{W}_{k}[n,s]:=\hat{S}_{k}[n,s]-kp$.  We refer to $\alp$ as the {\em centered hypergeometric distribution}.

Write $s=n(q-p)$ with $q=1-p$, and $N=np$ and $M=nq$. The support of $\alp$ is $\mathbb S:=\{a,a+1,\ldots ,b\}$ where $a=-[(kp)\wedge ((n-k)q)]$ and $b=(kq)\wedge ((n-k)p)$. Therefore $\alp$ has a Stein coefficient. For $x\in \mathbb S$
\bas
\alp(x) = \frac{\binom{N}{j}\binom{M}{k-j}}{\binom{n}{k}} \qquad \mb{ where }j=x+kp.
\eas
Then, 
\bas
\frac{\alp(x+1)}{\alp(x)}= \frac{(N-j)(k-j)}{(j+1)(M-k+j+1)}
\eas
which shows that the Markov chain with $\lam_{x}^{+}=(N-j)(k-j)=((n-k)p-x)(kq-x)$ and $\lam_{x}^{-}=j(M-k+j)=((n-k)q+x)(kp+x)$, is reversible for $\alp$. Since $\lam_{x}^{-}-\lam_{x}^{+}=nx$ by an easy calculation, we can write the rates as $\lam_{x}^{\pm}=T'(x)\mp \frac{n}{2}x$ where $T'(x)=\frac12(\lam_{x}^{+}+\lam_{x}^{-})$. Scale down all the rates by $n/2$ to get an Ehrenfest-like chain with rates $T(x)\mp x$ where
\bas
T(x)=\frac{\lam_{x}^{+}+\lam_{x}^{-}}{n} = 2pq\frac{k(n-k)}{n}+\frac{2x^{2}}{n}+(q-p)\frac{n-2k}{n}x.
\eas
Observe that if $k,p$ are fixed and $n\to \infty$, this converges to the Stein coefficient of the centered $\mb{Binomial}(k,p)$ distribution, as it should. When $p=\frac12$, the Stein coefficient takes the simpler form
\bas
T(x)=\frac{k(n-k)}{2n}+\frac{2x^{2}}{n}.
\eas

\eeg

\subsection{Stein coefficient under convolution} Let $X$ and $Y$ be independent random variables having Stein coefficients $T_{X}$ and $T_{Y}$ respectively and let $Z=X+Y$.   Let $f,g,h$ denote the probability mass functions of $X$, $Y$ and $Z$ respectively. Then $h(z)=\sum_{x}f(x)g(z-x)$ and by \eqref{eq:generalformulaforSteincoeff}, the Stein coefficient of $Z$ is 
\bas
T_{Z}(z) &= z+\frac{2}{h(z)}\sum_{w>z}h(w)w \\
&= z + \frac{2}{h(z)}\sum_{x,y:x+y>z}(x+y)f(x)g(y) \\
&= z + \frac{2}{h(z)}\l\{\sum_{x}f(x)\sum_{y: \ y>z-x}yg(y) \; + \; \sum_{y}g(y)\sum_{x: \ x>z-y}xf(x)\r\}.
\eas
Again using \eqref{eq:generalformulaforSteincoeff},  the two inner sums  are equal to $\frac12 g(z-x)(T_{Y}(z-x)-(z-x))$ and $\frac12 f(z-y)(T_{X}(z-y)-(z-y))$, respectively. Thus we get
\bas
T_{Z}(z) &= z +\frac{1}{h(z)}\sum_{(x,y):x+y=z}f(x)g(y)[T_{X}(x)-x+T_Y(y)-y] \\
&= \frac{1}{h(z)}\sum_{(x,y):x+y=z}f(x)g(y)[T_{X}(x)+T_{Y}(y)]
\eas
since $x+y=z$ for all $(x,y)$ in the sum and the sum of $f(x)g(y)$  precisely gives $h(z)$. We may write this as
\small
\begin{align}\label{eq:steincoeffofsumoftwo}
T_{Z}(z)=\E\l[T_{X}(X)+T_{Y}(Y)\Given X+Y=z\r].
\end{align}
\normalsize
More generally, if $X_{i}$ are independent  and $S_{n}=X_{1}+\ldots +X_{n}$, then the Stein coefficient $T_{n}$ of $S_{n}$ is given by
\small
\begin{align}\label{eq:steincoeffofsum}
T_{n}(s)=\E\l[\sum_{i=1}^{n}T_{X_{i}}(X_{i})\Given S_{n}=s\r].
\end{align}
\normalsize
\beg If $X_{i}$ are i.i.d. and take values $-p,q$ with probabilities $q,p$ respectively, then $T_{X_{i}}(x)=|x|$ for $x\in \{-p,q\}$.  From \ref{eq:steincoeffofsum}, the Stein coefficient of $S_{n}$ is $T_{n}(s)=\E[\sum_{i=1}^{n}|X_{i}|\Given S_{n}=s]$. But given $S_{n}=s$, then exactly $s+np$ of the $X_{i}$s are equal to $q$ and the remaining $nq-s$ are equal to $-p$.  Therefore, $\sum_{i=1}^{n}|X_{i}|=2pqn+(q-p)s$, reconfirming that $T_{n}(s)$ is what was found by direct calculation in Example~\ref{eg:binomialstein}.
\eeg

\subsection{Stein coefficient under scaling (and an additive perturbation)} Let $X$ be a zero mean random variable with support $\mathbb S=\{a,a+1,\ldots ,b\}$, mass function $f$ and Stein coefficient $T_{X}$. We wish to consider $2X$, but its support has gaps of length $2$. To get a random variable that has a Stein coefficient, we perturb it additively and consider $Y=2X+R$ where $R$ is a non-degenerate random variable having a Stein coefficient. What is the Stein coefficient of  $Y$? In Remark~\ref{rem:steinscalinggeneralR}, we record the result for general $R$, but for our purposes it suffices to take the simple case of $R$ taking the values $-1,0,1,$ with probabilities $\frac14,\frac12,\frac14$, respectively. 

For this choice of $R$, the support of  $Y$ is $\mathbb T=\{2a-1,2a,\ldots ,2b,2b+1\}$ and mass function $g(2x)=\frac12 f(x)$ and $g(2x-1)=\frac14 f(x)+\frac14 f(x-1)$. Hence, if $y=2x$, then by \eqref{eq:generalformulaforSteincoeff}
\bas
T_{Y}(y) &= y + \frac{2}{g(y)}\sum_{u>x}\l\{(2u-1)\frac{ f(u-1)+f(u)}{4} + 2u\frac{f(u)}{2} \r\}\\
&= y + \frac{4}{f(x)}\l\{ \frac14 (2x+1)f(x) + \sum_{u>x}f(u)\l[\frac{2u-1}{4}+\frac{2u+1}{4}+\frac{2u}{2}\r] \r\}\\
&=  y + (2x+1) + \frac{8}{f(x)}\sum_{u>x}uf(u).
\eas
The first two summands add up to $4x+1$ and the third summand is $4(T_{X}(x)-x)$, by \eqref{eq:generalformulaforSteincoeff}. Thus, $T_{Y}(y)=4T_{X}(x)+1$. 

Next suppose $y=2x-1$. Apply the defining formula \eqref{eq:Steinequationforbinomial} for $T_{Y}$ with $i=y$ and $i=y-1$ and average to get
\bas
T_{Y}(y)&=\frac{g(y+1)}{2g(y)}(T_{Y}(y+1)+y+1)+\frac{g(y-1)}{2g(y)}(T_{Y}(y-1)-(y-1)) \\
&=\frac{f(x)(4T_{X}(x)+y+2)+f(x-1)(4T_{X}(x-1)-y+2)}{f(x)+f(x-1)}
\eas
from the already worked out formulas for $T_{Y}$ at $y\pm 1$. In the above expression,  we may replace $f(x)$ by $T_{X}(x-1)-(x-1)$ and $f(x-1)$ by $T_{X}(x)+x$ (since the ratios are the same, by \eqref{eq:Steinequationforbinomial}), hence
\bas
T_{Y}(y) &= \frac{(T_{X}(x)-(x-1))(4T_{X}(x)+y+2)+(T_{X}(x)+x)(4T_{X}(x-1)-y+2)}{T_{X}(x)+T_{X}(x-1)+1} \\
&= \frac{B(4A-y)+A(4B+y)}{A+B} 
\eas
where $A=T_{X}(x)+x$ and $B=T_{X}(x-1)-(x-1)$. The numerator is $8AB+y(A-B)$, which can be written as $2(A+B)^{2}-2(A-B)^{2}+y(A-B)$. Therefore
\bas
T_{Y}(y) &= 2(A+B) - \frac{(A-B)(2(A-B)-y)}{A+B} \\
&=2(T_{X}(x)+T_{X}(x-1)+1) -R(y)
\eas
where 
\bas
R(y)=\frac{(y+[T_{X}(x)-T_{X}(x-1)])(y+2[T_{X}(x)-T_{X}(x-1)])}{T_{X}(x)+T_{X}(x-1)+1}.
\eas
In conclusion
\small
\begin{align}\label{eq:steincoeffunderscaling}
T_{Y}(y) =\begin{cases} 4T_{X}(x)+1 &\mb{ if }x=\frac12y\in \mathbb S, \\
2(T_{X}(x)+T_{X}(x-1)+1) -R(y) &\mb{ if }x=\frac12(y+1)\in \mathbb S. \end{cases}
\end{align}
\normalsize
Now we work out the two examples of centered binomial and hypergeometric distributions.
\beg\label{eg:binomialscaled} If $X$ is centered Binomial$(n,p)$, then we have seen in Example~\ref{eg:binomialstein} that $T_{X}(x)=2pqn+(q-p)x$ for $x\in \mathbb S=\{-np,-np+1,\ldots ,nq\}$. Therefore, $Y$ has Stein coefficient
\bas
T_{Y}(y)=\begin{cases} 8pqn+2(q-p)y+1 
 & \mb{ if }\frac12y\in \mathbb S,\\ 8pqn + 2(q-p)y+2-R_{n}(y)&\mb{ if }\frac12(y+1)\in \mathbb S.
\end{cases}
\eas
where $R_{n}(y)=\frac{(y+2q-2p)(y+q-p)}{4pqn+(q-p)y+1}$.  For the balanced case $p=\frac12$, 
\bas
T_{Y}(y)=\begin{cases} 2n+1 
 & \mb{ if }y\in \{-n,-n+2,\ldots ,n-2,n\},\\ 2n+ 2-\frac{y^{2}}{n+1} &\mb{ if }y\in \{-n-1,-n+1,\ldots ,n-1,n+1\}.
\end{cases}
\eas

\eeg

\beg \label{eg:hypergeoscaled} Let $X=\hat{W}_{k}[n,s]$ be the centered hypergeometric variable of Example~\ref{eg:hypergeometricstein} and let $Y=2X+R$. We have seen the formula for the Stein coefficient of $X$. It follows that if $y=2x-1$, then  $T_{X}(x-1)-T_{X}(x)=-\frac{2y}{n}-(q-p)\frac{n-2k}{n}$ and hence
\bas
R(y)=\frac{(y(1+\frac{2}{n})+(q-p)(1-\frac{2k}{n}))(y(1+\frac{4}{n})+2(q-p)(1-\frac{2k}{n})}{4pq\frac{k(n-k)}{n}+\frac{y^{2}+1}{n}+(q-p) (1-\frac{2k}{n})y+1}.
\eas
The expression looks complicated, but all that matters is that this is of smaller order than $y^{2}/k$ (as long as $k$ is away from $0$ and $n$).
\bas
T_{Y}(y)=\begin{cases} 8pq\frac{k(n-k)}{n}+\frac{2y^{2}}{n}+2(q-p) (1-\frac{2k}{n})y+1.
 & \mb{ if }\frac12 y\in \mathbb S, \\ 8pq\frac{k(n-k)}{n} +\frac{2y^{2}+2}{n} + 2 +2(q-p) (1-\frac{2k}{n})y- R(y)&\mb{ if }\frac12(y+1)\in \mathbb S.
\end{cases}
\eas
For the balanced case $p=\frac12$, we see that $R(y)=\frac{y^{2}(1+\frac{2}{n})(1+\frac{4}{n})}{\frac{k(n-k)}{n}+\frac{y^{2}+1}{n}+1}$ and
\bas
T_{Y}(y)=\begin{cases} 2\frac{k(n-k)}{n}+\frac{2y^{2}}{n}+1.
 & \mb{ if }y\in 2\mathbb S, \\ 2\frac{k(n-k)}{n} +\frac{2y^{2}+2}{n} + 2 - R(y)&\mb{ if }y\not\in 2\mathbb S.
\end{cases}
\eas
%

\eeg

\berk \label{rem:steinscalinggeneralR} For possible future use, we record the result for a general random variable $R$ that has mass function $r$ and Stein coefficient $T_{R}$. Then $Y=2X+R$ has Stein coefficient given by
\bas
T_{Y}(y) &= \frac{\summ_{2u+t=y}f(u)r(t)(2T_{X}(u)+T_{R}(t)) \; + \;  \summ_{2u+t=y\pm 1}f(u)r(t)(T_{X}(u)\pm u)}{\summ_{2u+t=y}f(u)r(t)}.
\eas
\eerk

\section{Proof of  Theorem~\ref{thm:bingausscouplingcauchy}}
Recall that $S_{n}$ is a sum of $n$  i.i.d. symmetric Bernoulli random variables. Then $\hat{S}_{n}:=\frac12 S_{n}$ has the centered $\mb{Bin}(n,\frac12)$ distribution. To couple $S_{n}$ and  $S_{4n}$ we use  the Markov chain coupling between  the Ehrenfest-like chains associated to $X=\hat{S}_{4n}$ and $Y=2\hat{S}_{n}+R=S_{n}+R$, where $R$ is as before (and independent of $S_{n}$).  If $\theta>0$, then 
\begin{align}\label{eq:bdforsns4nintermsofhattedones}
\E[e^{\theta|2S_{n}-S_{4n}|}] \le e^{2\theta}\E[e^{2\theta|2\hat{S}_{n}-\hat{S}_{4n}|}].
\end{align}

From the computation of Stein coefficients in Example~\ref{eg:binomialstein} and Example~\ref{eg:binomialscaled} (this is the case $p=\frac12$), we see that $Q(x,y)=|T_{Y}(y)-T_{X}(x)|$  is equal to $1$ if $y/2$ is in the support of $\hat{S}_{n}$ and $|2-\frac{y^{2}}{n+1}|$ otherwise. Thus, $Q\le \frac{Y^{2}}{n+1}+2$ and as $|Y|\le |S_{n}|+1$, we can write $Q\le \frac{2S_{n}^{2}+2}{n+1}+2\le \frac{2}{n}S_{n}^{2}+3$.

The Bernstein/Hoeffding inequality says that $\P\{S_{n}\ge t\}\le e^{-t^{2}/2n}$, which can be interpreted as saying that $(S_{n})_{+}^{2}$ and $(S_{n})_{-}^{2}$ are stochastically dominated by $2n\xi$, where $\xi$ is an exponential random variable with unit mean. Therefore, $\E[h(S_{n}^{2})]\le 2\E[h(2n\xi)]$ for any increasing $h:\R_{+}\mapsto \R$. In particular, we get $\E[h(Q)]\le 2\E[h(4\xi+3)]$.  In particular, for  $\alp<\frac{1}{4}$, we have 
\bas
\E[e^{\alp Q}] \; \le \; \frac{2e^{3\alp}}{1-4\alp} \; \le \; \frac{5}{1-4\alp}.
\eas
a uniform bound independent of $n$. 

Now, in the bound  \eqref{eq:bdfromMCcouplingforEhrenfest3} on $\E_{\gam}[e^{\theta|H|}]$ in the Markov coupling, one can take any $\theta$ small enough that $\frac{1}{2(1-\del)}e^{\theta}\theta^{2}<\frac{1}{4}$. For this, if  $2\theta^{2}e^{\theta}<1$, then we can always choose a $\del>0$ small enough. Then choose any $\mu$ to get a bound $\kappa$ which is independent of $n$. Thus, any $\theta_{0}$ with $2\theta_{0}^{2}e^{\theta_{0}}<1$ works. By \eqref{eq:bdforsns4nintermsofhattedones}, for the coupling of $S_{n}$ and $S_{4n}$, any $\theta_{0}$ with $8\theta_{0}^{2}e^{2\theta_{0}}<1$ suffices.
\hfill \qed

\section{Proof of the coupling Theorem~\ref{thm:hypergeogausscouplingcauchy}} 
As in the proof of Theorem~\ref{thm:bingausscouplingcauchy}, the variables $S_{k}[n,s]$ have spacing of $2$ in their support, hence we consider the modified variable $\hat{S}_{k}[n,s]=(S_{k}[n,s]+k)/2$ and $\hat{W}_{k}[n,s]=\hat{S}_{k}[n,s]-kp$ that were defined in Example~\ref{eg:hypergeometricstein}.  Let $\hat{W}_{1}=\hat{W}_{1}[n,0]$ and $\hat{W}_{2}=\hat{W}_{4k}[4n,0]$ and $\hat{W}=\hat{W}_{k}[n,s]$. These are half of $W_{1},W_{2},W$ that occur in the statement of Theorem~\ref{thm:hypergeogausscouplingcauchy}. Thus it suffices to prove the existence of a $\Theta>0$ such that there exist couplings satisfying
\small
\begin{align}
\E[e^{\theta|2\hat{W}_{1}-\hat{W}_{2}|}] &\le \frac32 \;\; \mb{ for all }\theta\le \Theta, \label{eq:reducedstatementhat1}  \\\
\E[e^{\theta|\hat{W}_{1}-\hat{W}|}] &\le e^{1+M\theta^{2}\frac{s^{2}}{n}} \;\; \mb{ for all }\theta\le \Theta. \label{eq:reducedstatementhat2} 
\end{align}
\normalsize

\bprf[Proof of \eqref{eq:reducedstatementhat1}] Let $V=\hat{W}_{1}$ and $X=\hat{W}_{2}$ and $Y=2V+R$ where $R$ takes values $-1,0,1$ with probabilities $\frac14,\frac12,\frac14$. In Example~\ref{eg:hypergeometricstein} and Example~\ref{eg:hypergeoscaled}, we have computed their  Stein coefficients
\bas
T_{X}(x) &= \frac{2k(n-k)}{n}+\frac{x^{2}}{2n}, \\ 
T_{Y}(y) &= \begin{cases} 2\frac{k(n-k)}{n}+\frac{2y^{2}}{n}+1
 & \mb{ if $y$ is even}, \\ 2\frac{k(n-k)}{n} +\frac{2y^{2}+2}{n} + 2 - \frac{y^{2}(1+\frac{2}{n})(1+\frac{4}{n})}{\frac{k(n-k)}{n}+\frac{y^{2}+1}{n}+1}&\mb{ if $y$ is odd}.
\end{cases}
\eas
The complicated looking term in $T_{Y}(y)$ is bounded as follows ($n\ge 2$):
\bas
\frac{y^{2}(1+\frac{2}{n})(1+\frac{4}{n})}{\frac{k(n-k)}{n}+\frac{y^{2}+1}{n}+1} \le \frac{27y^{2}}{n}
\eas
because $k(n-k)\ge \frac29 n^{2}$ (we just dropped the second and third terms in the denominator). Consequently, writing $y=2v+1$ (i.e., $v$ is the value of $V$) and using $y^{2}\le 8v^{2}+2$ and $k\le n$ we get
\bas
Q:=|T_{X}(x)-T_{Y}(y)| &\le C\l(\frac{x^{2}}{4k}+\frac{v^{2}}{k}+\frac{1}{k}\r)
\eas
for a constant $C$. Consider the bound \eqref{eq:bdfromMCcouplingforEhrenfest3} for the Markov coupling of $X$ and $Y$, with $H=|X-Y|$. Fix $\eps>0$ and set $\mu$ and $\del$ to be small positive numbers to get for some $M>0$, 
\bas
\E_{\gam}[e^{\theta|H|}] &\le (1+\eps) \E_{\gam}\l[e^{M\theta^{2}e^{\theta}Q}\r] \\
&\le (1+\eps)e^{\frac{1}{k}CM\theta^{2}e^{\theta}}\E\l[e^{2CM\theta^{2}e^{\theta} \frac{X^{2}}{4k}}\r]^{\frac12}\E\l[e^{2CM\theta^{2}e^{\theta} \frac{V^{2}}{k}}\r]^{\frac12}.
\eas
In the notation of Lemma~\ref{lem:boundsinbinomialandhypergeo}, $X=\frac12 W_{4k}$ and $V=\frac12 W_{k}$. Hence if $b=2 CMe^{\theta}\theta^{2}< \frac12$, part (1) of that Lemma shows that both expectations above are bounded by  $\frac{1}{\sqrt{1-4CM\theta^{2}e^{\theta}}}$. Therefore,
\ba
\E_{\gam}[e^{\theta|H|}] &\le (1+\eps)e^{\frac{1}{k}CM\theta^{2}e^{\theta}}\frac{1}{\sqrt{1-4CM\theta^{2}e^{\theta}}}.
\ea
As $|2\hat{W}_{1}-\hat{W}_{2}|\le |H|+1$, for $\theta>0$,
\ba
\E_{\gam}[e^{\theta|2\hat{W}_{1}-\hat{W}_{2}|}] &\le (1+\eps)e^{\frac{1}{k}CM\theta^{2}e^{\theta}}\frac{e^{\theta}}{\sqrt{1-4CM\theta^{2}e^{\theta}}}.
\ea
This can be made as close to $1$ as desired by taking $\theta\le \Theta$ for a small enough $\Theta>0$, and \eqref{eq:reducedstatementhat1} follows. 
\eprf

\bprf[Proof of \eqref{eq:reducedstatementhat2}] We want to find a good coupling between $\hat{W}_{1}=S_k[n,0]$ and $\hat{W}=S_k[n,s]-\frac{sk}{n}$. The Stein coefficients of these two variables were computed in Example~\ref{eg:hypergeometricstein} and Example~\ref{eg:hypergeoscaled} and are given by 
\bas
T_{1}(x) &= \frac12 \frac{k(n-k)}{n}+\frac{2x^{2}}{n}, \\
T(y) &= 2pq\frac{k(n-k)}{n}+\frac{2y^{2}}{n}+(q-p)\frac{n-2k}{n}y.
\eas
Therefore,
\bas
|T_{1}(x)-T(y)| &\le \frac12 (1-4pq)\frac{k(n-k)}{n} + \frac{2x^{2}}{n}+\frac{2y^{2}}{n} + |q-p|\ |1-\frac{2k}{n}| |y| \\
&\le \frac18 (p-q)^{2}n + \frac{2x^{2}}{n}+\frac{2y^{2}}{n} +\frac13 |q-p|\ |y| \\
&= \frac{s^{2}}{8n} +  \frac{2x^{2}}{n}+\frac{2y^{2}}{n} + \frac{|s|}{3\sqrt{n}}\frac{ |y| }{\sqrt{n}}.
\eas
Therefore, choosing $\del=1$ in \eqref{eq:bdfromMCcouplingforEhrenfest3} (it is here that we need this stronger form and \eqref{eq:bdfromMCcouplingforEhrenfest2} would not suffice), in the Markov coupling of $\hat{W}_{1}$ and $\hat{W}$, with $H=\hat{W}_{1}-\hat{W}$ (and recall that $\hat{W}_{1}=\frac12 W_{1}$ and $\hat{W}=\frac12 W$), 
\bas
\E[e^{\theta |H|}] &\le 4\E[e^{\theta^{2}|T_{1}(\hat{W}_{1})-T(\hat{W})|} \\
&\le 4e^{\theta^{2}\frac{s^{2}}{8n}}\E\l[e^{4\theta^{2}\frac{1}{n}\hat{W}_{1}^{2}} \r]^{\frac12} \E\l[e^{4\theta^{2} \frac{1}{n}\hat{W}^{2} + \frac{2\theta^{2} |s|}{3\sqrt{n}}\frac{1}{\sqrt{n}}|\hat{W}|} \r]^{\frac12} \\
&=4e^{\theta^{2}\frac{s^{2}}{8n}}\E\l[e^{\theta^{2}\frac{W_{1}^{2}}{n}} \r]^{\frac12} \E\l[e^{\theta^{2} \frac{W^{2}}{n} + \frac{\theta^{2} |s|}{3\sqrt{n}}\frac{|W|}{\sqrt{n}}} \r]^{\frac12} \\
&\le 4e^{\theta^{2}\frac{s^{2}}{8n}} \frac{1}{\sqrt[4]{1-8\theta^{2}}}\frac{2}{\sqrt[4]{1-8\theta^{2}}}e^{\frac{2\theta^{4}s^{2}}{9n(1-8\theta^{2})}}.
\eas
If $\theta\le \frac14$, then $1-8\theta^{2}\ge \frac12$ and the whole thing is  bounded by $12e^{\frac35 \frac{\theta^{2}s^{2}}{n}}$. Lastly $|W_{1}-W|=2|H|$, hence \eqref{eq:reducedstatementhat2} follows.
\eprf
%


\section{Some facts about the hypergeometric distribution}\label{sec:hypergeolemmas}
Consider a box that contains $n$ coupons, of which $np$ are labelled $+1$ and $nq$ are labelled $-1$. The sum of all the labels is $s:=n(p-q)$. Let $S_{k}'=X_{1}+\ldots +X_{k}$ where $X_{i}$s are drawn without replacement and let  $S_{k}=Y_{1}+\ldots +Y_{k}$ where $Y_{i}$s are drawn with replacement from the same box. We also write $W_{k}'=S_{k}'-k(p-q)=\sum_{i=1}^{k}(X_{i}-(p-q))$ and $W_{k}=S_{k}-k(p-q)=\sum_{i=1}^{k}(Y_{i}-(p-q))$ for the centered versions. Clearly, $W_{k}'$ has a centered hypergeometric distribution with parameters $(n,k,p)$ while $W_{k}$ has a centered Binomial distribution with parameters $(n,p)$.

We collect here certain results about these random variables, particularly bounds on expectations of certain functionals. While the Binomial is straightforward, the hypergeometric gets complicated if one tries a direct approach using the explicit mass function. Much difficulty can be avoided by using a wonderful result of Hoeffding~\cite{hoeffding} (this is the same paper where the famous Hoeffding inequality is proved, but this result is in a  somewhat less known Section~5 of the paper): 
\small
\begin{align}\label{eq:hoeffdingscomparison}
\E[f(S_{k}')]\le \E[f(S_{k})] \mb{ for any convex function }f:\R\mapsto \R.
\end{align}
\normalsize
Now we collect the results we need.
\begin{lemma}\label{lem:boundsinbinomialandhypergeo} 
\benu
\item For any real $\lam$ we have $\E[e^{\lam W_{k}'}] \; \le \;  \E[e^{\lam W_{k}}]  \; \le \; e^{\frac12 \lam^{2} k}$.
\item For any  real $a$ and $b<\frac12$, we have 
\bas
\E\l[e^{\frac{a}{\sqrt{k}}W_{k}'+\frac{b}{k}W_{k}'^{2}}\r] \; \le \;  \E\l[e^{\frac{a}{\sqrt{k}}W_{k}+\frac{b}{k}W_{k}^{2}}\r]  \; \le \; \frac{1}{\sqrt{1-2b}} e^{\frac{a^{2}}{2(1-2b)}}.
\eas
\item For any real $b<\frac12$, we have
\bas
\E\l[e^{\frac{b}{k}S_{k}'^{2}}\r] \; \le \;  \E\l[e^{\frac{b}{k}S_{k}^{2}}\r]  \; \le \; \frac{1}{\sqrt{1-2b}}e^{ \frac{b}{1-2b}\frac{k}{n}\frac{s^{2}}{n}}.
\eas
In particular, given any $\del>0$, there exists $b_{\del}>0$ and $c_{\del}<1$ such that for all $k\le (1-\del)n$ and any $b<b_{\del}$, these expectations are bounded by $\exp\{1+c_{\del}\frac{s^{2}}{n}\}$.
\eenu
\end{lemma}

\bprf The functions $x\mapsto e^{\lam (x-k(p-q))}$ and $x\mapsto e^{a(x-k(p-q))+b(x-k(p-q))^{2}}$ and $x\mapsto bx^{2}/k$ are all convex. By \eqref{eq:hoeffdingscomparison} this implies the first inequalities in all three statements of the lemma. It only remains to prove the second inequalities.
\benu
\item This is really a part of the proof of case of the famous Hoeffding's inequality, but let us quickly recap anyway. Since $Y_{i}-(p-q)$ takes the values $2q$ and $-2p$ with probabilities $p$ and $q$, respectively, we see that $\E[e^{\frac12\lam W_{k}}]=(pe^{\lam q}+qe^{-\lam p})^{k}$. Let $\phi(\lam)=\log (pe^{\lam q}+qe^{-\lam p})$ and observe that
\bas
\phi'(\lam)&= \frac{pq(e^{\lam q}-e^{-\lam p})}{pe^{\lam q}+qe^{-\lam p}} \mb{ vanishes at $\lam=0$, } \\
\phi''(\lam)&= \frac{pe^{\lam q}\times qe^{-\lam p}}{(pe^{\lam q}+qe^{-\lam p})^{2}}\le \frac14 \mb{ for any }\lam,
\eas
since $4ab\le (a+b)^{2}$ (here $a=qe^{-\lam p}$ and $b=pe^{\lam q}$). Thus,  $\phi(\lam)\le \frac18 \lam^{2}$ for all $\lam\in \R$, by the second order Taylor expansion. Replacing $\lam/2$ by $\lam$, we have arrived at $\E[e^{\lam W_{k}}] \le e^{\frac12 \lam^{2} k}$.

\item A useful and often used trick is to use the Parseval relation for characteristic functions to convert the quadratic in the exponent into a linear term by introducing a new independent Gaussian variable $Z$, independent of the other variables considered so far. As $\E[e^{\lam (Z+c)}]=e^{\lam c+\frac12 \lam^{2} }$, we can write
\bas
\E[e^{aW_{k}+bW_{k}^{2}}] &= \E[e^{aW_{k}+\sqrt{2b}W_{k}Z}] \\
&\le \E[e^{\frac12 k (a+\sqrt{2b}Z)^{2}}] \\
&=e^{\frac12 ka^{2}}\E[e^{ak\sqrt{2b}Z+bkZ^{2}}] \\
&=\frac{1}{\sqrt{1-2kb}}e^{\frac12 ka^{2}+\frac{k^{3}a^{2} b}{1-2kb}}
\eas
valid for any real $a$ and $b<\frac{1}{2k}$. Replace $a$ and $b$ by $a/\sqrt{k}$ and $b/k$ to get it in the form given in the statement of the Lemma. 
\item As $S_{k}=W_{k}+k(p-q)$, we see that
\bas
\E[e^{bS_{k}^{2}}] &= e^{bk^{2}(p-q)^{2}}\E[e^{bW_{k}^{2}+2bk(p-q)W_{k}}] \\
&\le \frac{1}{\sqrt{1-2kb}}\exp\l\{bk^{2}(p-q)^{2}+ \frac{4b^{2}k^{3}(p-q)^{2}}{2(1-2kb)}\r\} \\
&\le \frac{1}{\sqrt{1-2kb}}\exp\l\{ \frac{bk^{2}(p-q)^{2}}{1-2kb}\r\}.
\eas
Using $s=n(p-q)$ and replacing $b$ by $b/k$ we get
\bas
\E[e^{bS_{k}^{2}/k}] &=\frac{1}{\sqrt{1-2b}}\exp\l\{ \frac{b}{1-2b}\frac{k}{n}\frac{s^{2}}{n}\r\}.
\eas
This was the claimed inequality. Now it is clear that if $k\le n(1-\del)$, then choosing $b_{\del}<2\del$ ensures that the exponent is less than $c_{\del}s^{2}/n$ for $c_{\del}=(1-\del)/(1-2b_{\del})$. If $b<\frac14$, then $1/\sqrt{1-2b}$ is bounded by $\sqrt{2}<e$, hence the whole of it is bounded by $\exp\{1+c_{\del}\frac{s^{2}}{n}\}$. \hfill \qedhere
\eenu
\eprf

\newpage
\small

\section*{Appendix 1: From Tusnady type lemma to KMT theorem}
This derivation of KMT-EP from Tusnady's lemma is, up to changes of notation, copied from Pollard~\cite{pollard} (chapter 10) and  presented here for completeness. 

\smallskip
\noindent{\bf Dyadic intervals:} Let $\mathcal D=\cup_{p\ge 0}\mathcal D_{p}$, where $\mathcal D_{p}$ denotes the set of dyadic intervals of generation $p$, of the form $I=[k2^{-p},(k+1)2^{-p}]$ with $p\ge 0$ and $0\le k\le 2^{p}-1$.  Dyadic intervals have a natural rooted binary tree structure, with $I_{0}=[0,1]$ as the root and $I=[k2^{-p},(k+1)2^{-p}]$ having two children, $I'=[(2k)2^{-p-1},(2k+1)2^{-p-1}]$ and $I''=[(2k+1)2^{-p-1},(2k+2)2^{-p-1}]$. The ancestor of $I$ in the $q$th generation (for $0\le q\le p$) is denote $I^{q}$. Hence $I^{p}=I$ and $I^{0}=I_{0}$. The function $\psi_{I}:[0,1]\mapsto \R$ that is equal to $2^{p/2}$ on $I'$ and $-2^{p/2}$ on $I''$ and zero elsewhere, is called a Haar function. Together with the constant function $\one$, the Haar functions form an orthonormal basis of $L^{2}[0,1]$. The function $\phi_{I}(t):=\int_{0}^{t}\psi_{I}(s)ds$, vanishes outside $I$, takes the value $2^{-\frac12p -1}$ at the midpoint of $I$, and is linear on $I'$ and $I''$.

\smallskip
\noindent{\bf Series expansion of Brownian bridge:} Let $Z(I)$ be i.i.d. standard Gaussians. Then,   
$W_{0}^{(m)}(t):=\sum_{p\ge 0}^{m}\sum_{I\in \mathcal D_{p}}Z(I)\phi_{I}(t)$ converges uniformly (as $m\to \infty$) over $t\in [0,1]$ to the standard Brownian bridge $W_{0}$. If $t$ is a dyadic rational,  then $W_{0}^{(m)}(t)=W_{0}(t)$ for large $m$.

\smallskip
\noindent{\bf Constructing the uniform empirical process from $W_{0}$:} Fix $n$. Apply the coupling between Binomials and Gaussian assured by Lemma~\ref{lem:snandzcoupling} and construct $\{N(I)\suchthat I\in \mathcal D\}$ as follows.
 
Set $N(I_{0})=n$. Inductively, suppose $I$ is an interval for which $N(I)$ has been defined but $N(I')$ and $N(I'')$ are not yet defined. If $N(I)=0$, then set $N(I')=N(I'')=0$. If $N(I)\ge 1$, then couple $Z(I)$ with $\hat{N}(I)\eqd S_{N(I)}$ and set $N(I')=\frac12 (N(I)+\hat{N}(I))$ and $N(I'')=N(I)-N(I')$.  Inductively, it is easy to see that $N(I)\sim \mb{Bin}(n,|I|)$ where the length $|I|=2^{-p}$  if $I\in \mathcal D_{p}$.

From $N(I)$, $I\in \mathcal D$, we get $n$ i.i.d. uniform points in $[0,1]$, and $N$ is just the counting measure of these points. Let $\nu=\sqrt{n}(\frac1n N-\lam)$ (where $\lam$ denotes Lebesgue measure on $[0,1]$). The distribution function $G_{n}(t):=\nu[0,t]$, of the real measure $\nu$, is the uniform empirical process. Observe that $G_{n}'=\nu$, hence $\<G_{n}',\psi_{I}\>=\int_{0}^{1}\psi_{I}(t)d\nu(t)=\frac{1}{\sqrt{n}}2^{p/2}\hat{N}(I)$. Thus
$G_{n}(t)=\sum_{p\ge 0}\sum_{I\in \mathcal D_{p}}\frac{2^{p/2}}{\sqrt{n}}\hat{N}(I)\phi_{I}(t)$.

\noindent{\bf Closeness of the coefficients of $G_{n}^{(m)}$ and of $W_{0}^{(m)}$:}
The coupling in Lemma~\ref{lem:snandzcoupling} ensures that  $|\hat{N}(I)-Z(I)\sqrt{N(I)}|\le A(N(I),Z(I))$ and $|\hat{N}(I)|\le B(N(I),Z(I))$ where $A(n,t)=a(1+t^{2})$ and $B(n,t)=b(1+|t|\sqrt{n})$ (in fact the statement of Lemma~\ref{lem:snandzcoupling} is stronger).  For $I\in \mathcal D_{p}$,
\ba
& \big| \frac{2^{p/2}}{\sqrt{n}}\hat{N}(I)-Z(I) \big| \le \frac{2^{p/2}}{\sqrt{n}}|\hat{N}(I)-Z(I)\sqrt{N(I)}| \ + \ |Z(I)|\times \big|\sqrt{\frac{2^{p}}{n}N(I)}-1\big| \\
&\le \frac{2^{p/2}}{\sqrt{n}}A(N(I),Z(I)) \ + \ |Z(I)| \sum_{j=1}^{p}\Big| \sqrt{\frac{2^{j}}{n}N(I^{j})} - \sqrt{\frac{2^{j-1}}{n}N(I^{j-1})} \Big|
\ea
where in the second term we used the fact that $\frac{2^{0}}{n}N([0,1])=1$ to write a telescoping series. For $a,b>0$, we have $|\sqrt{a}-\sqrt{b}|\le |a-b|/\sqrt{b}$, because of which the $j$th summand in the second term can by bounded by  (assuming $N(I^{j-1})>0$, else the summand is zero)
\ba
&\frac{1}{\sqrt{\frac{2^{j-1}}{n}N(I^{j-1})}} \Big| \frac{2^{j}}{n}N(I^{j}) - \frac{2^{j-1}}{n}N(I^{j-1}) \Big| = \frac{2^{(j-1)/2}}{\sqrt{n}}\frac{|\hat{N}(I^{j-1})|}{\sqrt{N(I^{j-1})}} \\
&\le \frac{2^{(j-1)/2}}{\sqrt{n}}\frac{B(N(I^{j-1}),Z(I^{j-1}))}{\sqrt{N(I^{j-1})}}.
\ea
Plugging in the expressions for $A(n,t)$ and $B(n,t)$ and changing $j-1$ to $j$, 
\small
\begin{align}
&\big| \frac{2^{p/2}\hat{N}(I)}{\sqrt{n}}-Z(I) \big| \nonumber \\
&\le \frac{2^{p/2}a(1+|Z(I)|^{2})}{\sqrt{n}}  +   \frac{|Z(I)|}{\sqrt{n}}\sum_{j=0}^{p-1}2^{\frac{j}{2}}(\frac{b}{\sqrt{N(I^{j})}}+|Z(I^{j})|)\one_{N(I^{j})>0} \nonumber \\
&\le a\frac{2^{p/2}}{\sqrt{n}} + \frac{2^{p/2}}{\sqrt{n}}(|Z(I)|^{2}+\frac{b}{\sqrt{2}-1}|Z(I)|) \ + \  \frac{1}{2\sqrt{n}} \sum_{j=0}^{p-1}2^{\frac{j}{2}}(|Z(I)|^{2}+|Z(I^{j})|^{2}) \nonumber \\
&\le \frac{a'}{\sqrt{n}}2^{p/2} +  \frac{b'}{\sqrt{n}}\sum_{j=0}^{p}2^{j/2}|Z(I^{j})|^{2} \label{eq:boundfordiffofcoeffsGandW}
\end{align}
\small
for some constants $a',b'$.  In the second line we simply used $2xy\le x^{2}+y^{2}$ while to get to the last line, we first wrote $z^{2}+\frac{6}{\sqrt{2}-1}|z|\le 16(1+z^{2})$ and absorbed all the terms with $z^{2}$ into the last summand with $j=p$ (since the last term in the geometric series is of the same order as the sum). 

Fix $m$ and use this in the series expansion for $t\in T(m)$, where $T(m)$ is the set of end-points of the intervals in $\mathcal D_{m}$. Then, the series for $G_{n}$ and $W_{0}$ run up to $p=m$, and using \eqref{eq:boundfordiffofcoeffsGandW} and the fact that $|\phi_{I}(t)|\le 2^{-|I|/2}$, we get
\ba
|G_{n}(t)-W_{0}(t)| &\le \frac{a'}{\sqrt{n}}m +  \frac{b'}{\sqrt{n}}\sum_{p=0}^{m}\sum_{j=0}^{p}2^{(j-p)/2}|Z(I_{p,t}^{j})|^{2} \\
&\le \frac{a'}{\sqrt{n}}m  +\frac{4b'}{\sqrt{n}}\sum_{j=0}^{m}|Z(I_{j,t})|^{2}.
\ea
In the second line we used the fact that $I_{p,t}^{j}=I_{j,t}$ and interchanged the sums (the geometric series $\sum_{i}2^{-i/2}\le 4$). Hence, for $t\in T(m)$, we have $|G_{n}(t)-W_{0}(t)| = \frac{C}{\sqrt{n}}(m+S_{m}(t))$ 
where $S_{m}(t)=\sum_{p=0}^{m} |Z(I_{p,t})|^{2}$.  Pollard (for all references to \cite{pollard}, see section 10.7 of that book) shows that with $T(m)=\{k2^{-m}\suchthat 0\le k\le 2^{m}\}$, 
\small
\begin{align}\label{eq:bdSm}
\P\{\max_{t\in T(m)} S_{m}(t)\ge 10(m+x)\}\le 2e^{-m-x}. 
\end{align}
\small
This is easy to see: For a fixed $t\in T(m)$, the variable $S_{m}(t)$ has $\chi_{m+1}^{2}$ distribution. Hence $\E[e^{\lam S_{m}(t)}]=(1-2\lam)^{-\frac{m+1}{4}}$ for $\lam<\frac12$. Taking $\lam=\frac14$, we get $\P\{S_{m}(t)\ge 8(m+x)\}\le 2^{\frac14(m+1)-2(m+x)}$. The union bound gives something better than \eqref{eq:bdSm}.

Now let $\Del_{m}(f)=\max_{I\in \mathcal D_{m}}\max_{t,s\in I}|f(t)-f(s)|$ denote the maximum oscillation of a function $f$ within any interval in $ \mathcal D_{m}$. Then (see p. 254 of \cite{pollard}) using standard facts about Brownian bridge
\small
\begin{align}\label{eq:bdDelmW0}
\P\{\Del_{m}(W_{0})\ge \frac{m}{\sqrt{n}}\} \le 2^{m+1}e^{-\frac{2^{m}m^{2}}{2n}} \le 2^{m+1}e^{-\frac{m^{2}}{2}}.
\end{align}
Next, since $\sqrt{n}(G_{n}(t)-G_{n}(s))=(N_{n}(t)-N_{n}(s))-n(t-s)$ for $s<t$, 
\begin{align}
\P\{\Del_{m}(G_{n})\ge \frac{m+2}{\sqrt{n}}\} &\le \P\{\max_{I\in \mathcal D_{m}} N_{n}(I) \ge m+2-\frac{n}{2^{m}}\} \nonumber \\
&\le 2^{m}\P\{\mb{Bin}(n,2^{-m})\ge m\} \;\;  \le \frac{2^{m}}{m!}. \label{eq:bdDelmGn}
\end{align}
\small
as $\mb{Bin}(n,p)$ is stochastically dominated by $\mb{Pois}(np)$. Now, 
\ba
\max_{t\in [0,1]}|G_{n}(t)-W_{0}(t)| &\le \frac{C}{\sqrt{n}}\l\{  \max_{t\in T(m)}S_{m}(t) +\Del_{m}(G_{n})+\Del_{m}(W_{0})\r\}.
\ea
By \eqref{eq:bdSm}, \eqref{eq:bdDelmW0} and \eqref{eq:bdDelmGn}, one gets 
\ba
\P\l\{\max_{t\in [0,1]}|G_{n}(t)-W_{0}(t)|\ge C\frac{m+x}{\sqrt{n}} \r\} \le C'e^{-m-x}
\ea
which is the conclusion of KMT theorem, since $m=\log n+O(1)$. Usually it is written without the $m$ in the exponent, but as long as the constant $C$ is not specified in the left, it only makes sense to take $x\ge \log m$, which gives this term. It is neither stronger nor weaker to state it this way. 

\newpage

\section*{Appendix 2: Chatterjee's proof of KMT-RW  by induction}
For any probable value $s$ of $S_{n}$, let $\Sb[n,s]:=(S_{0}[n,s],\ldots ,S_{n}[n,s])$ (definition of $S_{k}[n,s]$ is in the paragraph preceding Lemma~\ref {lem:hypergeogausscoupling}. Let   $\Vb[n]:=(V_{0}[n],\ldots ,V_{n}[n])$ be a centered Gaussian vector with covariances $\E[V_{i}[n]V_{j}[n]]=\frac{i(n-j)}{n}$ for $0\le i\le j\le n$. When it is safe to do so without ambiguity, we shall drop the $n$ in the notation. 

\smallskip
\noindent{\bf Tools:} The following facts will be used.

\bllt By Lemma~\ref{lem:hypergeogausscoupling}, for any $n\ge 2$ and any $k\in [n/3,2n/3]$ and any $t$ in the support of $S_{n}$, there is a coupling such that for any $\theta\le \theta_{1}$, 
\begin{align}\label{eq:univhypergeocoupling}
\E\l[e^{\theta|S_{k}[n,t]-\frac{k}{n}t-V_{k}[n]|}\r] \le e^{1+M\theta^{2}\frac{t^{2}}{n}}. 
\end{align}

\bllt By part (3) of Lemma~\ref{lem:boundsinbinomialandhypergeo}, there exists $\alp_{0}>0$ and $\gam<1$ such that for any $n\ge 1$ and any $k\le \frac23 n$, any probable value $t$ of $S_{n}$, and any $\alp\le \alp_{0}$, 
\begin{align}\label{eq:expsquareofhypergeo}
\E\l[e^{\frac{\alp}{k} S_{k}[n,t]^{2}}\r] \le e^{1+\gam \alp\frac{t^{2}}{n}}.
\end{align}

Fix constants $A,B,\lam_{0}$ satisfying: $A\ge \frac{1+\log 2}{\log \frac32}$ and  $B\ge \frac{2M}{1-\gam}$ and  $\lam_{0}\le \frac{\theta_{1}}{2}\wedge \sqrt{\frac{\alp_{0}}{2B}}$ where $M,\gam,\theta_{1}$ are those occurring in \eqref{eq:univhypergeocoupling} and \eqref{eq:expsquareofhypergeo}.

\para{Induction hypothesis} For any $n\ge 1$ and any probable value $t$ of $S_{n}$, there is a coupling of $\Sb[n,t]$ with $\Vb[n]$ such that for  any $\lam\le \lam_{0}$, 
\begin{align}\label{eq:inductionhypothesis}
\E\l[\exp\l\{\lam \max\limits_{1\le i\le n}|S_{i}[n,t]-\frac{it}{n}-V_{i}[n]|\r\}\r] \le e^{A\log n+B\lam^{2}\frac{t^{2}}{n}}.
\end{align}

\smallskip
\noindent{\bf Base case:} If $A$ is large enough, then the statement is obvious for any fixed $n$. We just need to watch out that the induction step goes through without difficulty if $A$ is increased.

\smallskip
\noindent{\bf Induction step:} Assume that the conclusion \eqref{eq:inductionhypothesis} holds for any $m<n$ in place of $n$ (and any probable value $t$ of $S_{m}$).  Fix $n\ge 6$ and a probable value $t$ of $S_{n}$. Choose $k=\lfloor n/2\rfloor$ so that $\frac{n}{3}\le k\le \frac{2n}{3}$ as $n\ge 6$. 

\setlength{\leftmargini}{0pt}
\begin{enumerate}
\item Construct  $(s,v)$ having the same marginals as $(S_{k}[n,t],V_{k}[n])$ and coupled so that \eqref{eq:univhypergeocoupling} holds. That is, with $R=|s-\frac{k}{n}t-v|$, we have $\E[\exp\{\theta R\}]\le \exp\{1+M\theta^{2}\frac{t^{2}}{n}\}$ for all $\theta\le \theta_{1}$.
\item Conditional on $(s,v)$, construct independent pairs $(\Sb',\Ub')$ and $(\Sb'',U'')$ such that $\Sb'\eqd \Sb[k,s]$, $\Ub'\eqd\Vb[k]$, $\Sb''\eqd \Sb[n-k,t-s]$ and $\Ub''\eqd \Vb[n-k]$ so that \eqref{eq:inductionhypothesis} holds for both pairs. That is, writing
\ba
T'= \max\limits_{1\le i\le k}|S'_{i}-\frac{is}{k}-U'_{i}|\} \;\; \mb{ and }\;\; T''=\max\limits_{1\le j\le n-k}|S''_{j}-\frac{j(t-s)}{n-k}-U''_{j}|,
\ea
for all $\lam\le \lam_{0}$ we have
\ba
&\E\l[e^{\lam T'}\Given s,v \r] \le e^{A\log k+B\lam^{2}\frac{s^{2}}{k}}  \;\; \mb{ and } \;\; &\E\l[e^{\lam T''}\Given s,v\r] \le e^{A\log (n-k)+B\lam^{2}\frac{(t-s)^{2}}{n-k}}.
\ea
\item Define $\Sb=(S_{0},\ldots ,S_{n})$ and $\Vb=(V_{0},\ldots ,V_{n})$ by setting
\ba
S_{i}=\begin{cases} S'_{i} &\mb{ if }i\le k, \\ s+S''_{i-k} & \mb{ if }k\le i\le n, \end{cases}  \;\; \mb{ and } \;\;
V_{i}=\begin{cases} U'_{i}+\frac{i}{k}v &\mb{ if }i\le k, \\ U''_{i-k}+\frac{n-i}{n-k}v & \mb{ if }k\le i\le n. \end{cases}
\ea
There is no ambiguity at $i=k$, as $S'_{k}=s$ and $S''_{0}=U'_{k}=U''_{0}=0$. From the construction, it is clear that $\Sb\eqd \Sb[n,t]$ and $\Vb\eqd \Vb[n]$.
\item Now we observe that
\ba
S_{i}-\frac{it}{n}-V_{i} &= \begin{cases}
(S'_{i}-\frac{is}{k}-U'_{i}) +\frac{i}{k}(s-\frac{kt}{n}-v) &\mb{}i\le k, \\
(S''_{i-k}-\frac{(i-k)(t-s)}{n-k}-U''_{i-k}) +\frac{n-i}{n-k}(s-\frac{kt}{n}-v) &\mb{}k\le i\le n.
\end{cases}
\ea
Consequently, if $T:=\max_{i\le n}|S(i)-\frac{i}{n}a-V(i)|$ then $T\le ( T'\vee T'')+R$ and hence $\E[e^{\lam T}] \le \E[e^{\lam (T'+R)}]+\E[e^{\lam (T''+R)}]$. Now, for any $p\in (0,1)$ (to be chosen depending on $\gam$), we can write
\ba
\E[e^{\lam (T'+R)}] &= \E\l[e^{\lam R}\E[e^{\lam T'} \Given s,v]\r] \\
&\le e^{A\log k} \ \E\l[e^{\lam R}e^{B\lam^{2}\frac{s^{2}}{k}}\r] \\
&\le e^{A\log k} \ \E\l[e^{2\lam R}\r]^{\frac12}\E\l[e^{2B\lam^{2}\frac{s^{2}}{k}}\r]^{\frac12}
\ea
by H\"{o}lder's inequality. Since $2\lam\le \theta_{1}$ and $2B\lam^{2}\le \alp_{0}$, 
\ba
\E[e^{\lam (T'+R)}] &\le e^{A\log k} \ e^{\frac12(1+4M\lam^{2}\frac{t^{2}}{n})} \ e^{\frac12(1+2B\lam^{2}\gam  \frac{t^{2}}{n})} \\
&= e^{A\log k} \ e^{1+\frac{\lam^{2}t^{2}}{n}[2M+\gam B]} \\
&\le e^{-A\log\frac{n}{k}+1} \ e^{A\log n+B\lam^{2}\frac{t^{2}}{n}}
\ea
provided $2M\le (1-\gam)B$. By almost identical reasoning, we also get
\ba
\E[e^{\lam (T''+R)}] &\le  e^{-A\log\frac{n}{n-k}+1} \ e^{A\log n+B\lam^{2}\frac{t^{2}}{n}}.
\ea
 Therefore,
\ba
\E[e^{\lam T}] &\le e^{A\log n+B\lam^{2}\frac{t^{2}}{n}}\l(e^{-A\log\frac{n}{k}+1}+e^{-A\log\frac{n}{n-k}+1} \r) \\
&\le e^{A\log n+B\lam^{2}\frac{t^{2}}{n}}
\ea
since the condition $k\le \frac23 n$ ensures that second factor is bounded by $2e^{1-A\log\frac32}<1$.
\eenu
This completes the induction step.

\smallskip
\noindent{\bf The KMT-RW theorem for symmetric Bernoulli steps} From \eqref{eq:inductionhypothesis}, it is easy to deduce the KMT-RW theorem for symmetric Bernoullis (Chatterjee proves something stronger, keeping the same Brownian motion as $n$ varies, but we forego that strengthening now). For this, let $W$ be a standard Brownian motion, and recall that we may write $W(t)=W_{0}(t)+tZ$ where $Z$ is a standard Gaussian independent of the Brownian bridge $W_{0}$. Further, $\Vb[n]$ has the same distribution as $W_{0}$ sampled at times $0,\frac1n,\ldots ,\frac{n-1}{n},1$. Now the idea is clear.
\setlength{\leftmargini}{0pt}
\begin{enumerate}
\item Invoke Lemma~\ref{lem:bingausscoupling} and couple $S_{n}$ and $Z$ so that $\E[e^{\theta|S_{n}-z\sqrt{n}|}]\le \kappa_{0}$ for any $\theta\le \theta_{0}$.
\item Conditional on $Z=z$, the distribution of $W$ is $W_{0}(t)+tz$ and conditional on $S_{n}=s$, the distribution of $(S_{0},\ldots ,S_{n})$ is the same as $\Sb[n,s]$. Now we can couple $\Sb$ with the random vector $\Vb=(W_{0}(0),W_{0}(1/n),\ldots, W_{0}(1))$ so that \eqref {eq:inductionhypothesis} holds.
\item We observe that
\ba
\max_{0\le k\le n}|S_{k}-W(k)|\le \max_{0\le k\le n}|S_{k}-\frac{k}{n}S_{n}-W_{0}(k/n)| \; + \; |S_{n}-W(1)|.
\ea
Hence $\E[\exp\{\lam \max_{0\le k\le n}|S_{k}-W(k)|\}]$ is bounded by  
\ba
&\E[\exp\{\lam \max_{0\le k\le n}|S_{k}-\frac{k}{n}S_{n}-W_{0}(k/n)|\}] + \E[\exp\{ \lam |S_{n}-W(1)|\}] \\
&\le e^{A\log n}\E[e^{B\lam^{2}S_{n}^{2}/n}] + \kappa_{0}
\ea
if $\lam\le \lam_{0}\wedge \theta_{0}$. By Lemma~\ref{lem:boundsinbinomialandhypergeo}, if $\lam_{0}$ is sufficiently small, then this whole quantity is bounded by $\kappa e^{A\log n}$ for some constants $\kappa$ and $A$ that do not depend on $n$. By Markov's inequality, we get KMT-RW. 
\end{enumerate}

%

\bigskip

\end{document}